\documentclass{birkmult}
\usepackage{amssymb}
\usepackage[OT2,T2A]{fontenc}
\newtheorem{theorem}{Theorem}[section]

\newtheorem{corollary}[theorem]{Corollary}
\newtheorem{lemma}[theorem]{Lemma}
\newtheorem{proposition}[theorem]{Proposition}

 \theoremstyle{definition}
 
 \theoremstyle{remark}

 \numberwithin{equation}{section}

\begin{document}
\title{On maximal proper subgroups of field automorphism groups} 
\author{M.Rovinsky}
\address{Independent University of Moscow \\
119002 Moscow \\ B.Vlasievsky Per. 11\\ 
and Institute for Information \\ 
Transmission Problems \\ 
of Russian Academy of Sciences}
\email{\tt marat@mccme.ru} 
\thanks{The author gratefully acknowledges the support of the Alexander 
von Humboldt-Stiftung and of Pierre Deligne 2004 
Balzan prize in mathematics. The author was supported in part by 
RFBR grants 06-01-72550-CNRS-L{\fontencoding{OT2}\selectfont \_a}, 
07-01-92211-CNRS-L{\fontencoding{OT2}\selectfont \_a}} 


\date{January 1, 2004}

\begin{abstract} Let $G$ be the automorphism group of an extension 
$F|k$ of algebraically closed fields of characteristic zero of 
transcendence degree $n$, $1\le n\le\infty$. In this paper we 
\begin{itemize} \item construct some maximal closed non-open 
subgroups $G_v$, and some (all, in the case of countable 
transcendence degree) maximal open proper subgroups of $G$; 
\item describe, in the case of countable transcendence 
degree, the automorphism subgroups over the intermediate 
subfields (a question of Krull, \cite[\S4, question 3b)]{krull}); 
\item construct, in the case $n=\infty$, a fully 
faithful subfunctor $(-)_v$ of the forgetful functor 
from the category ${\mathcal S}m_G$ of smooth representations 
of $G$ to the category of smooth representations of $G_v$; 
\item construct, using the functors $(-)_v$, a subfunctor $\Gamma$ 
of the identity functor on ${\mathcal S}m_G$, coincident (via the 
forgetful functor) with the functor $\Gamma$ on the category of 
{\sl admissible} semilinear representations of $G$ constructed in 
\cite{adm} in the case $n=\infty$ and $k=\overline{{\mathbb Q}}$. 
\end{itemize} The study of open subgroups is motivated by 
the study of (the stabilizers of the) smooth representations 
undertaken in \cite{repr,adm}. The functor $\Gamma$ is an analogue 
of the global sections functor on the category of sheaves on a 
smooth proper algebraic variety. Another result is that `interesting' 
semilinear representations are `globally generated'. \end{abstract}
\maketitle 

Consider a field extension $F|k$ and a group $H$ of its automorphisms. 

Given a presheaf (i.e. a contravariant functor) ${\mathcal F}$ on a 
category of smooth $k$-varieties (whose morphisms include all smooth ones) 
one can associate to it the `generic fibre' ${\mathcal F}(F):=
\lim\limits_{_U\longrightarrow\phantom{_U}}{\mathcal F}(U)$, where 
$U$ runs over the smooth irreducible $k$-varieties with function fields 
embedded into $F$. It is endowed with a natural $H$-action, and should 
be considered as one of typical objects studied in birational geometry. 

In most general terms, we are going to study relations between 
presheaves ${\mathcal F}$ and their `generic fibres' ${\mathcal F}(F)$. 

Following \cite{jac,pss,sh,ih} (and generalizing the case \cite{krull-gal} 
of algebraic extension), consider $H$ as a topological group with the 
base of open subgroups given by the stabilizers of finite subsets 
of $F$. Then $H$ becomes a totally disconnected Hausdorff group. 

Any element of the `generic fibre' ${\mathcal F}(F)$ is the image 
of an element of ${\mathcal F}(U)$ for some $U$, and therefore, is 
fixed by the subgroup of $H$ leaving the function field of $U$ fixed. 
This means that the $H$-set ${\mathcal F}(F)$ is {\sl smooth}, i.e. 
its stabilizers are open. 

Evidently, one cannot get much of a representation of a `small' 
group, so we fix the following setting. 

Let $F|k$ be an extension of countable or finite transcendence 
degree $n$, $1\le n\le\infty$, of algebraically closed fields of 
characteristic zero, and $G=G_{F|k}$ be its automorphism group. 

Then for any intermediate subfield $k\subseteq K\subseteq F$ 
the topology on $G_{F|K}$ coincides with the restriction of the 
topology on $G$. $G$ is locally compact if and only if $n<\infty$. 

In \S\ref{mono-morfizm} we describe the open subgroups of $G$ in 
terms of certain `primitive' open subgroups of $G_{L|k}$ for 
algebraically closed extensions $L|k$ of finite transcendence degree. 

In \S\ref{decomp-podgr} we study discrete valuations of $F$ trivial 
on $k$ and the corresponding decomposition groups. This is applied in 
\S\ref{assoc-funktory} to the `globalization', i.e. to constructing 
some presheaves on some categories of smooth $k$-varieties, of smooth 
$G$-modules. Because of the links with some types of motives, the emphasis 
is made on the subcategory ${\mathcal I}_G$ of smooth $G$-modules (see 
\S\ref{stabilizatory}) and related categories of semilinear representations. 

In \S\ref{globalizaciq} to a discrete valuation $v$ of $F$ trivial on $k$ a 
fully faithful functor $(-)_v$ is associated. It is not clear yet, whether 
$(-)_v$ is exact. Its exactness is related to the following geometric 
question, cf. Lemma \ref{tochnostt}. Do there exist \begin{itemize} 
\item an integer $N\ge 1$, \item an irreducible variety $X$ over $k$, 
\item a collection of $N$ surjective maps $f_j:X\longrightarrow Y_j$, 
$1\le j\le N$, with $\dim Y_j<\dim X$ \end{itemize}
such that the induced homomorphism of the 0-cycles 
$Z_0(X)\longrightarrow\bigoplus_{j=1}^NZ_0(Y_j)$ is injective?

\section{Structure of $G$} 
The classical morphism $\beta:\{\mbox{subfields $F$ over $k$}\}
\hookrightarrow\{\mbox{closed subgroups of $G$}\}$, given by 
$K\mapsto G_{F|K}$, is injective, inverts the inclusions, transforms 
the compositum of subfields to the intersection of subgroups, and 
respects the units: $k\mapsto G$. The image of $\beta$ is stable 
under the passages to sup-/sub- groups with compact quotients; 
$\beta$ identifies the subfields over which $F$ is algebraic 
with the compact subgroups of $G$ (\cite{jac,pss,sh,ih}). 

In particular, the proper subgroups in the image of 
$\beta$ are the compact subgroups in the case $n=1$. 

The map $H\mapsto F^H$, left inverse of $\beta$, 
inverts the order, but it is not a monoid morphism. 

Let $A\Pi$ be the set of algebraically closed subfields of $F$ of 
finite transcendence degree over $k$. There is a morphism of partially 
ordered commutative associative unitary monoids (transforming 
the intersection of subgroups to the algebraic closure of the 
compositum of subfields) $\alpha:\{\mbox{open subgroups of $G$}\}
\longrightarrow\!\!\!\!\rightarrow A\Pi$, uniquely determined 
by the condition $G_{F|\alpha(U)}\subseteq U$ and the 
transcendence degree of $\alpha(U)$ over $k$ is minimal. 
Details are in Proposition \ref{morfizm-mono}. 

\begin{theorem}[\cite{repr}] \label{top-simpl} \begin{enumerate} 
\item The subgroup $G^{\circ}$ of $G$, generated by the compact 
subgroups, is open and topologically simple, if $n<\infty$. 
If $n=\infty$ then $G^{\circ}$ is dense in $G$. 
\item \label{prosto-besk} Any closed normal proper subgroup of 
$G$ is trivial, if $n=\infty$. \end{enumerate} \end{theorem}

If $n<\infty$, the left $G$-action on the one-dimensional 
oriented ${\mathbb Q}$-vector space of right-invariant measures 
on $G$ gives rise to a surjective homomorphism, the modulus, 
$\chi:G\longrightarrow{\mathbb Q}^{\times}_+$, which is 
trivial on $G^{\circ}$. However, I know nothing about 
the discrete group $\ker\chi/G^{\circ}$. If it is trivial for 
$n=1$ then it is trivial for arbitrary $n<\infty$. If $n=1$ 
a presentation of $G/G^{\circ}$ as a quotient of a quite 
structured discrete group is given in Lemma \ref{sur-fil}. 

\vspace{4mm}

Then one can characterize the image of $\beta$ 
in the case $n=\infty$ in the following 4 steps. 
\begin{enumerate} \item \label{normaliz-alg-zamk}
The normalizers $G_{\{F,L\}|k}$ of $G_{F|L}$ in $G$ for all 
$L\in A\Pi\smallsetminus\{k\}$ are the maximal open proper 
subgroups of $G$. (This follows from Proposition \ref{morfizm-mono}.) 
\item \label{ochen-zamk}
The subgroups $G_{F|L}$ of $G$ for all $L\in A\Pi\smallsetminus\{k\}$ 
are minimal closed non-trivial normal subgroups in $G_{\{F,L\}|k}$ 
from (\ref{normaliz-alg-zamk}). (This follows from Theorem 
\ref{top-simpl} (\ref{prosto-besk}). Namely, let $H\subset G_{\{F,L\}|k}$ 
be a closed non-trivial normal subgroup and $1\neq\sigma\in H$. 
Let us show that $H$ intersects $G_{F|L}$ non-trivially. One has 
$\sigma\tau\sigma^{-1}\tau^{-1}\in H\cap G_{F|L}$ for any $\tau\in G_{F|L}$. 
By \cite[Corollary 2.3]{repr} there exists $\tau\in G_{F|L}$ such that 
$\sigma\tau\sigma^{-1}\tau^{-1}\neq 1$, which means $H\cap G_{F|L}\neq\{1\}$. 
Now Theorem \ref{top-simpl} (\ref{prosto-besk}) implies that $H$ 
contains $G_{F|L}$.) 
\item \label{otkr-kon-tipa}
The subgroups $G_{F|L}$ of $G$ for all non-trivial extensions 
$L\subset F$ of $k$ of finite type are the open subgroups containing 
normal co-compact subgroups of type $G_{F|\overline{L}}$ from 
(\ref{ochen-zamk}). (This follows from the above Galois theory 
for $\overline{L}|k$.) 
\item 
The proper subgroups in the image of $\beta$ are intersections 
of subgroups from (\ref{otkr-kon-tipa}). 
\end{enumerate}

{\it Remark.} The subgroups $G_{F|L}$ of $G$ for all extensions 
$L\subset F$ of $k$ of finite type and transcendence degree one 
are the subgroups $G_{F|L}$ from (\ref{otkr-kon-tipa}) 
with the only maximal proper subgroup of $G$ containing them. 

The above procedure can be modified as follows. 
\begin{enumerate} \item \label{normaliz-alg-zamk-str}
The subgroups $G_{\{F,\overline{k(x)}\}|k}$ of $G$ for all 
$x\in F\smallsetminus k$ are the open subgroups $U$ of $G$ 
such that $\#(U\backslash G/U)=2$. 
\item \label{ochen-zamk-str}
The subgroups $G_{F|\overline{k(x)}}$ of $G$ for all 
$x\in F\smallsetminus k$ are minimal closed non-trivial normal 
subgroups in $G_{\{F,\overline{k(x)}\}|k}$ 
from (\ref{normaliz-alg-zamk-str}). 
\item \label{otkr-kon-tipa-str}
The subgroups $G_{F|k(x)}$ of $G$ for all $x\in F\smallsetminus k$ 
are the open subgroups containing normal co-compact subgroups of 
type $G_{F|\overline{k(x)}}$ from (\ref{ochen-zamk-str}) 
and such that $N_GU/U\cong{\rm PGL}_2k$. 
\item 
The proper subgroups in the image of $\beta$ are intersections 
of subgroups from (\ref{otkr-kon-tipa-str}). 
\end{enumerate}
This is based on the following two observations. (i) The condition 
$\#(U\backslash G/U)=2$ means that $G=U\sqcup U\sigma U$ for some 
(and therefore, any) $\sigma\in G\smallsetminus U$, and thus, 
$U$ is maximal among proper subgroups of $G$. 
(ii) For any $F',F''\in A\Pi\smallsetminus\{k\}$ the set 
$G_{\{F,F'\}|k}\backslash G/G_{\{F,F''\}|k}$ is a disjoint 
union of the sets $S_j$, 
$0\le j\le\min({\rm tr.deg}(F'|k),{\rm tr.deg}(F''|k))$, 
where $$S_j=G_{\{F,F'\}|k}\backslash\{K\subset F~|~\mbox{$K\cong F''$ 
over $k$},~{\rm tr.deg}(KF'|F')={\rm tr.deg}(F''|k)-j\}.$$

\section{Maximal open subgroups} \label{mono-morfizm}
In this section we assume that ${\rm tr.deg}(F|k)=n\le\infty$.
\begin{lemma} \label{red-tuda} If $\overline{L_1},\overline{L_2}$ 
are proper subfields of $F$ and $1<{\rm tr.deg}(\overline{L_1}|
\overline{L_1}\cap\overline{L_2})<\infty$ then the group $H$, generated 
by $G_{F|\overline{L_1}}$ and $G_{F|\overline{L_2}}$, contains 
$G_{F|\overline{L}}$, where $\overline{L}$ is a subfield of $F$, not 
containing $L_1$, and such that $\overline{L_1}\cap\overline{L}$ 
contains $\overline{L_1}\cap\overline{L_2}$ as a proper subfield. 
\end{lemma} 
\begin{proof} Set $k:=\overline{L_1}\cap\overline{L_2}$. 
Clearly, $H$ contains $G_{F|\sigma(\overline{L_1})}$
for any $\sigma\in G_{F|\overline{L_2}}$. As 
$\sigma(\overline{L_1})$ is algebraically closed and 
$\overline{L_1}\cap\sigma(\overline{L_1})\supseteq k$, 
it remains to choose $\sigma\in G_{F|\overline{L_2}}$ such 
that $\sigma^{-1}(\overline{L_1})\not\subseteq\overline{L_1}$ 
and $\overline{L_1}\cap\sigma(\overline{L_1})\neq k$, 
and to set $L:=\sigma(\overline{L_1})$. 

Let $x\in\overline{L_1}\smallsetminus k$ and $y\in\overline{L_1}
\smallsetminus\overline{k(x)}$. We need $\sigma\in G_{F|\overline{L_2}}$ 
such that $\sigma x\in\overline{L_1}$ and $\sigma^{-1}y\not\in
\overline{L_1}$. If $y$ is not in $\overline{L_2(x)}$ then such 
$\sigma\in G_{F|\overline{L_2}}$, clearly, exists. (The 
$G_{F|\overline{L_2(x)}}$-orbit of $y$ spans $F$.) 

Suppose now that the element $y$ is forced to be algebraic over $L_2(x)$, 
i.e., ${\rm tr.deg}(L_1L_2|L_2)=1$. Then $x$ and $y$ are related 
by a polynomial $P(x,y)=0$ with coefficients in $\overline{L_2}$. Let 
us enumerate the non-zero coefficients of $P$ by a set $\{0,\dots,N\}$. 
We shall assume that the polynomial $P(X,Y)=\sum_{s=0}^Np_sX^{i_s}Y^{j_s}$ 
is irreducible over $\overline{L_2}$, and $p_0=1$. Suppose that 
$\sigma^{-1}(\overline{L_1})\subseteq\overline{L_1}$ (which is equivalent 
to $\sigma(\overline{L_1})=\overline{L_1}$, since 
${\rm tr.deg}(\overline{L_1}|k)<\infty$) for any 
$\sigma\in G_{F|\overline{L_2}}$ such that $\sigma x\in\overline{L_1}$. 
Then $-(dP)(\sigma x,\sigma y)=P_I(\sigma x,\sigma y)d(\sigma x)
+P_{II}(\sigma x,\sigma y)d(\sigma y)\in
F\otimes_{L_1}\Omega^1_{L_1|k}\subseteq\Omega^1_{F|k}$, which can be 
rewritten as $\sum_{s=1}^N\sigma x^{i_s}\sigma y^{j_s}dp_s
\in F\otimes_{L_1}\Omega^1_{L_1|k}$. It remains to show that the 
whole space $F^N$ is the linear envelope of the vectors 
$(\sigma x^{i_1}\sigma y^{j_1},\dots,\sigma x^{i_N}\sigma y^{j_N})$ for 
$\sigma\in G_{F|\overline{L_2}}$ and $\sigma x\in\overline{L_1}$ (since 
then $dp_s\in F\otimes_{L_1}\Omega^1_{L_1|k}$ for all $s$, which means that 
$p_s\in\overline{L_1}$, i.e., $p_s\in \overline{L_1}\cap\overline{L_2}=:k$).

Otherwise, if such vectors belong to a hyperplane: 
$\sum_{s=1}^N\lambda_s\sigma x^{i_s}\sigma y^{j_s}=0$ for all 
$\sigma\in G_{F|\overline{L_2}}$ as above, then we may assume 
that $\lambda_1=1$ and the number $M$ of non-zero $\lambda_s$ is 
minimal: $\sum_{s=1}^M\lambda_s\sigma x^{i_s}\sigma y^{j_s}=0$.
Subtracting from this equality the image of the equality 
$\sum_{s=1}^M\lambda_sx^{i_s}y^{j_s}=0$ under the action of 
$\sigma$, we get that $\lambda_s\in F^{G_{F|\overline{L_2}}\cap 
G_{\{F,\overline{L_1}\}|k}}$. As $F^{G_{F|\overline{L_2}}\cap 
G_{\{F,\overline{L_1}\}|k}}\subseteq\bigcap
_{z\in\overline{L_1}\smallsetminus k}\overline{L_2}(z)\subseteq
\bigcap_{b=1}^{\infty}\overline{L_2}(x^b)=\overline{L_2}$, 
we got thus another polynomial relation between $x$ and $y$ over 
$\overline{L_2}$, contradicting the minimality of the polynomial $P$. 
\end{proof}

\begin{corollary} \label{perv-sled}
If ${\rm tr.deg}(L_1|L_1\cap L_2)<\infty$, $L_1,L_2$ are proper 
algebraically closed subfields in $F$ and $L_1$ is not contained 
in $L_2$ then $H\supseteq G_{F|L}$ for some algebraically closed 
subfield $L$ such that ${\rm tr.deg}(L_1|L_1\cap L)=1$. \end{corollary}
\begin{proof} The set of algebraically closed subfields of $L$,
not containing $L_1$ and such that $H\supseteq G_{F|L}$, is non-empty, 
since contains $L_2$. Let ${\rm tr.deg}(L_1|L_1\cap L)$ be minimal.
According to Lemma \ref{red-tuda}, ${\rm tr.deg}(L_1|L_1\cap L)=1$. 
\end{proof}

\begin{lemma} \label{red-cor-odi} If $L_1,L_2$ are algebraically 
closed subfields in $F$ and ${\rm tr.deg}(L_1|L_1\cap L_2)=1$ 
then the subgroup $H$ generated by $G_{F|L_1}$ and $G_{F|L_2}$ 
is dense in $G_{F|L_1\cap L_2}$. \end{lemma}
\begin{proof} This is a version of Lemma 2.16 from \cite{repr}. We may 
assume that $L_2$ is a proper subfield of $F$ not contained in $L_1$. 
Thus, $L_1$ is an algebraic closure of $k(x)$ for some 
$x\in F\smallsetminus L_2$, where $k:=L_1\cap L_2$. Then for any 
$z\in F\smallsetminus L_2$, $y\in F\smallsetminus L_1$ there exist 
$\sigma\in G_{F|L_2}$, $\tau\in G_{F|L_1}$ such that $\sigma x=z$, 
$\tau x=y$. Then for any $y\in F\smallsetminus k$ 
there exists $\sigma\in G_{F|L_1}G_{F|L_2}$ such that $\sigma x=y$. 
It follows that $H$ is a normal subgroup of $G$. If $n=\infty$ then 
the topological group $G$ is simple (Theorem 2.9, \cite{repr}), 
i.e., $H$ is dense in $G$. If $n<\infty$ then, according to the 
same Theorem 2.9, \cite{repr}, the topological group $G^{\circ}$ 
is also simple, and thus, the closure of $H$ contains $G^{\circ}$. 
It is known (Lemma 2.15, \cite{repr}) that $G=G_{F|L_1}G^{\circ}$, 
i.e., the closure of $H$ again coincides with $G$. \end{proof} 

\begin{corollary} \label{sled-2} If $L_1,L_2$ are proper algebraically 
closed subfields of $F$, $L_1$ is not contained in $L_2$ and 
${\rm tr.deg}(L_1|L_1\cap L_2)<\infty$ then the closure of $H$ contains 
$G_{F|L}$ for some proper algebraically closed subfield $L$ in $L_1$. 
\end{corollary}
\begin{proof} According to Corollary \ref{perv-sled}, there exists 
an algebraically closed subfield $L'$ in $F$ such that ${\rm tr.deg}
(L_1|L_1\cap L')=1$ and $H\supseteq G_{F|L'}$. It follows from Lemma 
\ref{red-cor-odi} that the subgroup generated by $G_{F|L_1}$ and 
$G_{F|L'}$ is dense in $G_{F|L}$, where $L:=L_1\cap L'$. We deduce 
from this that the closure of $H$ contains some $G_{F|L}$ of 
the desired type. \end{proof} 

\begin{proposition} \label{2.14} The subgroup $H=\langle G_{F|L_1},
G_{F|L_2}\rangle$ is dense in $G_{F|L_1\cap L_2}$ for any subfields 
$L_1$ and $L_2$ in $F$ such that $\overline{L_1}\cap\overline{L_2}$ is 
algebraic over $L_1\cap L_2$ and ${\rm tr.deg}(L_1|L_1\cap L_2)<\infty$. 
\end{proposition}
\begin{proof} The usual Galois theory reduces the problem to the case of 
algebraically closed $L_1$ and $L_2$. Next, we may assume that $L_2$ does 
not contain $L_1$, and $L_1\neq F$. Set $k:=L_1\cap L_2$. Let 
$L=\overline{L}\supseteq k$ be such that $H\supseteq G_{F|L}$ and let 
${\rm tr.deg}(L|k)(\le{\rm tr.deg}(L_1|k))$ be minimal. If 
$L\not\subseteq L_j$ ($j=1$ or $j=2$) then, according to Corollary 
\ref{sled-2}, (${\rm tr.deg}(L|L\cap L_j)\le{\rm tr.deg}(L|k)<\infty$), 
the closure of the subgroup generated by $G_{F|L}$ and 
$G_{F|L_j}$, (and thus, the closure of $H$ as well) contains 
$G_{F|L'}$, where $L'$ is a proper algebraically closed subfield in $L$, 
contradicting the minimality of ${\rm tr.deg}(L|k)$. Thus, $L=k$. \end{proof}

\begin{proposition} \label{morfizm-mono} 
There is a morphism of commutative associative unital monoids 
inverting inclusions (transforming the intersection of subgroups 
to the algebraic closure of the compositum of subfields, and the 
identity $G$ to the identity $k$) $$\left\{\begin{array}{c}
\mbox{{\rm open}}\\ \mbox{{\rm subgroups of} $G$}\end{array}
\right\}\stackrel{\alpha}{\longrightarrow\!\!\!\!\!\rightarrow}\left
\{\begin{array}{c}\mbox{{\rm algebraically closed subfields of} $F$}\\
\mbox{{\rm of finite transcendence degree over} $k$}\end{array}
\right\}=:A\Pi,$$ uniquely determined by the following equivalent 
conditions: \begin{itemize} \item $G_{F|\alpha(H)}$ is a normal 
subgroup of $H$ 
and, if possible, $\alpha(H)\neq F$; 
\item $G_{F|\alpha(H)}\subseteq H$ and ${\rm tr.deg}(\alpha(H)|k)$ 
is minimal. \end{itemize} In particular, for any $L\in A\Pi$ 
distinct from $F$ and $k$, the normalizer $G_{\{F,L\}|k}$ of 
$G_{F|L}$ is maximal among proper subgroups of $G$. 

If $n=\infty$ then any proper open subgroup $H$ of $G$ is contained 
in a maximal proper subgroup of $G$, and any maximal proper open 
subgroup of $G$ is of type $G_{\{F,L\}|k}$ for some $L\in A\Pi$, 
$L\neq k$. Besides that, $\alpha(H)=\alpha(N_GH)$. \end{proposition}
\begin{proof} If a subgroup $H$ of $G$ is open then there exists 
$L\in A\Pi$ such that $G_{F|L}\subseteq H$. Let $L$ be of minimal 
possible transcendence degree. For any $\sigma\in H$ the group 
$H$ contains the closure of the subgroup generated by $G_{F|L}$ 
and $G_{F|\sigma(L)}$. Then, by Proposition \ref{2.14}, $H$ contains 
$G_{F|L\cap\sigma(L)}$, and therefore, $\sigma(L)\supseteq L$. Thus, 
$H\subseteq G_{\{F,L\}|k}$. In other words, $G_{F|L}\triangleleft H$.

Let us check the maximality of $G_{\{F,L\}|k}$ for $L\neq k,F$. 
If $G_{\{F,L\}|k}\subseteq U$ then $G_{F|K}\subseteq U\subseteq 
G_{\{F,K\}|k}$ for some $K=\overline{K}$ such that ${\rm tr.deg}
(K|k)\le{\rm tr.deg}(L|k)$, since $G_{F|L}\subseteq U$. The inclusion 
$G_{\{F,L\}|k}\subseteq G_{\{F,K\}|k}$ takes place only if either 
$K=L$, or $K=k$. In the second case $U=G$. 

If, moreover, $G_{F|L'}\subseteq H$ for another $L'\in A\Pi$ then 
either $L\subseteq L'$, or $L=F$. If, in addition, $H\subseteq 
G_{\{F,L'\}|k}$ then either $L=L'$, or $L=F$, or $L'=F$. This 
shows that $\alpha$ is well-defined. This also implies that 
$G_{F|\alpha(H)}\triangleleft N_GH$, and thus, 
$\alpha(H)=\alpha(N_GH)$, if $n=\infty$. 

Let us check that $\alpha$ is a morphism. If $G_{F|\overline{K}}
\subseteq U\subseteq G_{\{F,\overline{K}\}|k}$ and 
$G_{F|\overline{L}}\subseteq V\subseteq G_{\{F,\overline{L}\}|k}$ 
then $G_{F|\overline{K}~\overline{L}}\subseteq U\cap V\subseteq 
G_{\{F,\overline{K},\overline{L}\}|k}:=G_{\{F,\overline{K}\}|k}
\bigcap G_{\{F,\overline{L}\}|k}$. As 
$G_{F|\overline{K}~\overline{L}}\supseteq G_{F|\overline{KL}}$ 
and $G_{\{F,\overline{K},\overline{L}\}|k}\subseteq 
G_{\{F,\overline{KL}\}|k}$, it remains to show that 
$\alpha(U\cap V)=F$ in the case $\overline{KL}=F$. 

If $\alpha(U)=F$ then evidently $\alpha(U\cap V)=F$. If 
$\overline{K}$ and $\overline{L}$ are distinct from $F$ 
then $\alpha(G_{\{F,\overline{K},\overline{L}\}|k})=
\overline{KL}$, i.e. again $\alpha(U\cap V)=F$. \end{proof} 

\vspace{4mm}

{\sc Remarks.} 1. If $n<\infty$ and $H\subset G$ is contained in 
neither subgroup of type $G_{\{F,L\}|k}$ for $L\in A\Pi\smallsetminus\{k,F\}$ 
then $F$ is algebraic over the subfield generated over $k$ by the 
$H$-orbit of $x$ for any $x\in F\smallsetminus k$. 

2. If ${\rm tr.deg}(L|k)={\rm tr.deg}(F|L)=\infty$ then 
$G_{\{F,\overline{L}\}|k}$ is maximal among the proper closed subgroups 
of $G$, i.e., the subgroup $H$ generated by $G_{\{F,\overline{L}\}|k}$
and by any $\sigma\in G$ such that $\sigma(\overline{L})\neq\overline{L}$ 
is dense in $G$. {\it Question.} Can one replace the condition 
`${\rm tr.deg}(F|L)=\infty$' by the condition `$F\neq\overline{L}$'? 

\begin{proof} Replacing if necessary $\sigma$ by $\sigma^{-1}$, we may 
assume that $\sigma(\overline{L})$ does not contain $\overline{L}$. 
According to \cite[Proposition 2.14]{repr}, the closure of $H$ contains 
$G_{F|\overline{L}\cap\sigma(\overline{L})}$, since this is the closure 
of the subgroup generated by $G_{F|\overline{L}}$ and 
$G_{F|\sigma(\overline{L})}$. Set $L'=\overline{L}\cap
\sigma(\overline{L})$. 

Then there exists an element $\tau\in G_{\{F,\overline{L}\}|k}$ such that 
$L'\cap\tau(L')=k$. Namely, choose a transcendence base $t_1,t_2,t_3,\dots$ 
of $\overline{L}|k$ such that $L'\subseteq\overline{k(t_2,t_3,\dots)}$. Let 
$\tau\in G_{\{F,\overline{L}\}|k}$ be such that $\tau t_j=t_1^{j-1}+t_j$ 
for any $j\ge 1$. Then any element of $L'\cap\tau(L')$ has the following 
form $F(t_1+t_2,t_1^2+t_3,\dots,t_1^{M-1}+t_M)=G(t_2,t_3,\dots,t_M)$ for 
some algebraic functions $F$ and $G$, and some $M\ge 2$. The equality 
of the exterior differentials in $\Omega^1_{F|k}$ gives us 
$\sum_{j=1}^{M-1}F_jd(t_1^j+t_{j+1})=\sum_{j=1}^{M-1}G_jdt_{j+1}$, where 
$F_j=\partial F/\partial X_j$ and $G_j=\partial G/\partial X_j$, or 
$(\sum_{j=1}^{M-1}jt_1^{j-1}F_j)dt_1=\sum_{j=1}^{M-1}(G_j-F_j)dt_{j+1}$,
which is equivalent to $F_j=G_j$ for all $1\le j<M$ and 
$\sum_{j=1}^{M-1}jt_1^{j-1}G_j=0$, i.e., $F_j=G_j=0$ for all 
$1\le j<M$. This means that $L'\cap\tau(L')=k$. \end{proof} 

3. Suppose that finitely generated subfields $K_1,\dots,K_N$ of $F|k$ 
are in general position. Then the common stabilizer of their algebraic 
closures $G_{\{F,\overline{K_1},\dots,\overline{K_N}\}|k}$ is contained in 
precisely $2^N-1$ maximal proper open subgroups of $G$: $G_{\{F,K_S\}|k}$ 
for all non-empty subsets $S\subseteq\{1,\dots,N\}$, where $K_S$ is 
the algebraic closure of the compositum of $K_j$ for $j\in S$. 
Indeed, $F,k$ and $K_S$ for all $S$ are the only algebraically closed 
subfields preserved by $G_{\{F,\overline{K_1},\dots,\overline{K_N}\}|k}$. 

4. For any totally disconnected group $H$ and a subgroup 
$U\subset H$ the functors $H^0(U,-)$ and $\lim\limits_{U\overrightarrow
{\subset V\subset}H}H^0(V,-)$ on the category of smooth $H$-sets (or 
modules, etc.) coincide. Here $V$ runs over the open subgroups of $H$. 

In particular, the following conditions are equivalent 
\begin{enumerate} \item \label{invarianty-mod} $H^0(U_1,-)=H^0(U_2,-)$ 
on the category of smooth representations of $H$, 
\item \label{invarianty-mnozh} 
$H^0(U_1,-)=H^0(U_2,-)$ on the category of smooth $H$-sets, 
\item \label{otkrytye-podgruppy} any open subgroup of $H$ containing $U_1$ 
contains $U_2$ and vice versa. \end{enumerate} 
\begin{proof} Implications (\ref{otkrytye-podgruppy})$\Rightarrow$%
(\ref{invarianty-mnozh})$\Rightarrow$(\ref{invarianty-mod}) are trivial. 
(\ref{invarianty-mod})$\Rightarrow$(\ref{otkrytye-podgruppy}). Assuming 
the contrary, let an open subgroup $V$ contains $U_1$ but not $U_2$. 
Then $[1]\in H^0(U_1,{\mathbb Q}[H/V])$, but 
$[1]\not\in H^0(U_2,{\mathbb Q}[H/V])$. \end{proof} 

Then one can characterize the union of proper open subgroups of $G$ as follows. 
\begin{corollary} \label{harakteriz-dopolnen} If $n=\infty$ then 
the union of proper open subgroups of $G$ is a proper dense subset of $G$. 
The following properties of an element $\sigma\in G$ are equivalent. 
\begin{enumerate} \item \label{nelezhit} 
$\sigma$ does not belong to the union of the proper open subgroups of $G$, 
\item \label{G-inv-mnozh} $W^{\langle\sigma\rangle}=W^G$ 
for any smooth (i.e., with open stabilizers) $G$-set $W$, 
\item \label{G-inv-preds} $W^{\langle\sigma\rangle}=W^G$ for any 
smooth representation $W$ of $G$, \item \label{G-inv-q-form} 
there are non-zero $\sigma$-invariant finite-dimensional 
$F$-vector subspaces in $\Omega^q_{F|k}$ for neither $q\ge 1$, 
\item \label{G-inv-1-form} there are no non-zero $\sigma$-invariant 
finite-dimensional $F$-vector subspaces in $\Omega^1_{F|k}$. 
\end{enumerate}
If $\sigma$ has these equivalent properties then for any smooth $G$-group 
$M$ restriction to $\langle\sigma\rangle$ of any non-trivial smooth 
$G$-torsor under $M$ is non-trivial, i.e., the sequence of pointed 
sets $\{\ast\}\longrightarrow H^1_{Sm}(G,M)\longrightarrow 
H^1(\langle\sigma\rangle,M)$ is exact.\footnote{but 
$H^1_{Sm}(G,M)\hookrightarrow H^1(\langle\sigma\rangle,M)$ is not 
bijective, as shows the example of trivial $G$-group $M\neq\{1\}$.} 
As one can it imagine, the groups ${\rm Hom}_{\langle\sigma\rangle}$ are 
in general much bigger than the groups ${\rm Hom}_G$. \end{corollary}
\begin{proof} For instance, if $\{x_i~|~i\in{\mathbb Z}\}$ 
is a transcendence base of $F$ over $k$ and $\sigma\in G$ is such 
that $\sigma x_i=x_{i+1}$ for any $i\in{\mathbb Z}$ then $G$ 
is the only open subgroup containing $\sigma$. 

To show the density, note that for any $\tau\in G$ and any $L\subset F$ 
of finite transcendence degree over $k$ there is $\sigma_L\in 
G_{\{F,\overline{L\tau(L)}\}|k}$ such that $\sigma_L|_L=\tau|_L$. 

(\ref{nelezhit}), (\ref{G-inv-mnozh}), (\ref{G-inv-preds}) 
are equivalent by the preceding remark applied to $U_1=G$ 
and $U_2=\langle\sigma\rangle$. 

(\ref{G-inv-mnozh})$\Rightarrow$(\ref{G-inv-q-form}). Consider 
the set of non-zero finite-dimensional $F$-vector subspaces in 
$\Omega^q_{F|k}$. Clearly, it is smooth, and it has no elements 
fixed by the whole group $G$. Then any element fixing a point of 
our set belongs to a proper open subgroup. 

The implication (\ref{G-inv-q-form})$\Rightarrow$(\ref{G-inv-1-form}) is 
trivial. The implication (\ref{G-inv-1-form})$\Rightarrow$(\ref{G-inv-mnozh}) 
follows from Proposition \ref{morfizm-mono}. 

Let $K=k(x_i~|~i\in{\mathbb Z})$ and $K^+=k(x_i~|~i\ge 0)$. 
The Haar measure on $G_{F|K}$ induces an embedding 
${\rm Hom}_{\langle\sigma\rangle}(K,W)\hookrightarrow
{\rm Hom}_{\langle\sigma\rangle}(F,W)$ via the normalized trace 
surjection $F\longrightarrow K$. Let $M=\{\frac{P}{Q}\in K^+
~|~-1\le\deg_{x_0}P<\deg_{x_0}Q,~P,Q\in\sigma(K^+)[x_0]\}$. 
Then the the natural map $k\oplus{\mathbb Q}[\langle\sigma\rangle]
\otimes\left(M\oplus x_0\sigma(K^+)[x_0]\right)\longrightarrow K$ 
is an isomorphism of $\langle\sigma\rangle$-modules, and thus, 
${\rm Hom}_{\langle\sigma\rangle}(K,W)={\rm Hom}(k,W^G)\oplus
{\rm Hom}(M\oplus x_0\sigma(K^+)[x_0],W)$. \end{proof} 

\vspace{4mm}

{\sc Remarks.} 1. Even if there is a $q$-dimensional $\sigma$-invariant 
subspace in $\Omega^1_{F|k}$ there does not always exist a 
$\sigma$-invariant subfield in $F$ of transcendence degree $q$ over $k$. 
{\sc Example.} Let $x$, $y$, $z_i$, $i\in{\mathbb Z}$, form a transcendence 
base of $F|k$, and $\sigma x=y$, $\sigma y=x^2$, 
$\sigma z_i=z_{i+1}$. Then $\sigma\omega=-\sqrt{2}\omega$, where 
$\omega=\sqrt{2}\frac{dx}{x}-\frac{dy}{y}$. If there is a $\sigma$-invariant 
subfield in $F$ of transcendence degree $1$ over $k$, and $\varphi$ is its 
non-constant element then consider $d\varphi=\varphi_Idx+\varphi_{II}dy$. 
Then $d\varphi\wedge d\sigma\varphi=0$, so 
$\varphi_I\cdot\sigma\varphi_I=2x\varphi_{II}\cdot\sigma\varphi_{II}$. 
Set $\psi=\varphi_I/\varphi_{II}$. Then $\psi\cdot\sigma\psi=2x$. 
As the equation $\alpha\cdot\sigma\alpha=1$ has only constant solutions, 
$\alpha=\pm 1$ (since this implies that $\alpha=\sigma^2\alpha$), 
and therefore, $\psi=\pm\sqrt{2}y/x$. As the analytic solutions of 
$x\varphi_I(x,y)-\sqrt{2}y\varphi_{II}(x,y)=0$ are of type 
$\varphi(x,y)=\xi(x^{\sqrt{2}}y)$, we get $\varphi\in\overline{k}$. \qed 

2. If $n=\infty$ then any countable free group $H=\ast_{j\in S}{\mathbb Z}$ 
can be embedded into $G$ in such a way that its intersection with any 
proper open subgroup in $G$ is trivial. Namely, choose a transcendence 
base of $F|k$, and enumerate it by the elements of $H$: $\{x_h~|~h\in H\}$. 
Define an action of the generators $\{h_j~|~j\in S\}$ of $H$ on the 
transcendence base by $h_jx_h=x_{h_jh}$. Clearly, this action 
extends, though not uniquely, to $F$. 

\vspace{4mm}

{\it Questions.} 1. The preimage of any subgroup of a prime index 
in ${\mathbb Q}^{\times}_+$ under the modulus character gives, 
if $n<\infty$, an example of a maximal open proper subgroup, not 
encounted by Proposition \ref{morfizm-mono}. Any compact subset of 
$G$ is contained in infinitely many subgroups of this type. 
Are there any other maximal proper open subgroups? 

2. Do there exist closed subgroups not contained in maximal proper ones?

3. Let $A_{\sigma}:=F\langle\sigma,\sigma^{-1}\rangle$ be the algebra 
of endomorphisms of the additive group $F$ generated by $F$ and by 
$\sigma^{\pm 1}$ for some $\sigma\in G$. Clearly, $A_{\sigma}$ is a 
Euclidean simple central $F^{\langle\sigma\rangle}$-algebra, cf. \cite{O}. 
The set of $\sigma$-invariant algebraically closed subfields in $F|k$ 
injects into the set of $A_{\sigma}$-submodules in $\Omega^1_{F|k}$ by 
$L\mapsto F\otimes_L\Omega^1_{L|k}$. 
Suppose that $n=\infty$ and $\sigma$ does not belong to the union of 
the proper open subgroups of $G$. In particular, $\Omega^1_{F|k}$ is a 
torsion-free $A_{\sigma}$-module of at most countable rank. In a 
standard manner one checks that the finitely generated torsion-free 
$A_{\sigma}$-modules are free. 

An example of the $A_{\sigma}$-module $\Omega^1_{F|k}$ of rank 1, 
which is not free, is given by 
$F=\overline{k(x_i,y_j~|~i\in{\mathbb Z},~j\in{\mathbb N})}$, where $x_i,y_j$ 
are algebraically independent, we set $y_0=x_0$ and $\sigma x_i=x_{i+1}$, 
$\sigma y_j=y_{j-1}+y_j$. In particular, there is a strictly increasing 
sequence of $\sigma$-invariant algebraically closed subfields in $F|k$: 
$\overline{k(\sigma^{{\mathbb Z}}x_0)}\subset
\overline{k(\sigma^{{\mathbb Z}}y_1)}\subset
\overline{k(\sigma^{{\mathbb Z}}y_2)}\subset\dots$

The rank of the $A_{\sigma}$-module $\Omega^1_{F|k}$ is an invariant 
of the conjugacy class of $\sigma$. What are the others? 

In the case $F=\overline{k(x_i~|~i\in{\mathbb Z})}$, where $x_i$ are 
algebraically independent and $\sigma x_i=x_{i+1}$, one has 
$A_{\sigma}\stackrel{\sim}{\longrightarrow}\Omega^1_{F|k}$, 
$\alpha\mapsto\alpha dx_0$, so the set of $A_{\sigma}$-submodules in 
$\Omega^1_{F|k}$ is in bijection with the set of left ideals in 
$A_{\sigma}$, i.e., with the set of monic (non-commutative) polynomials 
in $\sigma$ with non-zero constant term. E.g., the polynomial $\sigma+1$ 
corresponds to $\overline{k(x_i+x_{i+1}~|~i\in{\mathbb Z})}\neq F$. 
The $\sigma^{{\mathbb Z}}$-orbit of an element $y\in F$ form a 
transcendence base of $F|k$ if and only if there exists $i\in{\mathbb Z}$ 
such that $y\in\overline{k(x_i)}\smallsetminus k$. 

\subsection{Stabilizers of `homotopy invariant' representations} 
\label{stabilizatory}
Denote by ${\mathcal I}_G$ the full subcategory in ${\mathcal S}m_G$ 
consisting of representations $W$ such that $W^{G_{F|M}}=W^{G_{F|M'}}$ 
for any extension $M$ of $k$ in $F$ and any purely transcendental 
extension $M'$ of $M$ in $F$. Denote by ${\mathcal I}={\mathcal I}_{/k}=
{\mathcal I}_{F|k}:{\mathcal S}m_G\longrightarrow{\mathcal I}_G$ the left 
adjoint to the inclusion functor 
${\mathcal I}_G\hookrightarrow{\mathcal S}m_G$, and set 
$C_L:={\mathcal I}_{F|k}{\mathbb Q}[\{L\stackrel{/k}{\hookrightarrow}F\}]$ 
for any extension $L$ of $k$ of finite type, cf. \cite{repr}, \S6. 

The stabilizers of any object of ${\mathcal I}_G$ are of type 
$G_{\{F,F'\}|k}\times_{G_{F'|k}}H$, where $F'\in A\Pi$ (is an algebraically 
closed extension of $k$ in $F$ of finite transcendence degree) and $H$ is 
an open subgroup of $G_{F'|k}$ such that if $G_{F'|L}\subseteq H$ then 
$\overline{L}=F'$. (This follows from Proposition \ref{morfizm-mono}, since 
if ${\rm Stab}_w=G_{\{F,F_1,\dots,F_s\}|k}\times_{G_{F_s|k}}H$, where 
$H\subseteq G_{F_s|k}$ is open and if $G_{F_s|L}\subseteq H$ then 
$\overline{L}=F_s$; $F_1\supseteq\dots\supseteq F_s$ is a flag in $A\Pi$, 
and $S$ is the union of transcendence bases of $F_1|F_2,\dots,F_{s-1}|F_s$ 
then ${\rm Stab}_w\supseteq G_{F|L(S)}$, and therefore, 
${\rm Stab}_w\supseteq G_{F|L}$, i.e. $F_1=F_s$.) 

One might wonder, whether the objects of ${\mathcal I}_G$ (or of 
${\mathcal A}dm_G$) are determined by their stabilizers. The following 
example suggests that the morphisms are not determined by the stabilizers, 
if one allows irreducible objects. 

{\sc Example.} There exist linear maps, respecting the stabilizers which 
are not $G$-equivariant. Namely, let $A$ be a commutative algebraic $k$-group, 
$\alpha$ be a ${\mathbb Q}$-linear homomorphism $A(F)/A(k)\longrightarrow 
A(k)_{{\mathbb Q}}$, and $\Lambda\subseteq A(k)_{{\mathbb Q}}$ be a 
${\mathbb Q}$-vector subspace. Then 
$id+\alpha:W:=A(F)_{{\mathbb Q}}/\Lambda\longrightarrow W$ preserves 
the stabilizers, but it is $G$-equivariant only if $\alpha=0$. 

\section{Valuation subgroups} \label{decomp-podgr}
Let ${\mathcal O}_v$ be a valuation ring in $F$, $\mathfrak{m}_v
={\mathcal O}_v\smallsetminus{\mathcal O}_v^{\times}$ be the maximal ideal, 
and $\kappa(v)$ be the residue field. If $k\subseteq{\mathcal O}_v$, 
fix a subfield $k\subseteq F'\subseteq{\mathcal O}_v$ identified 
with $\kappa(v)$ by the reduction modulo $\mathfrak{m}_v$. In this 
case $\kappa(v)$ is of characteristic zero (and algebraically closed). 

Set $G_v:=\{\sigma\in G~|~\sigma({\mathcal O}_v)=
{\mathcal O}_v\}$. This is a closed subgroup in $G$. 

The ${\mathbb Q}$-vector space $\Gamma:=F^{\times}/{\mathcal O}
_v^{\times}$ is totally ordered: $v(x)\ge v(y)$ if and only if 
$xy^{-1}\in{\mathcal O}_v$, where $v:F^{\times}\longrightarrow
\Gamma$ is the natural projection. The rank of $v$ is 
$r:=\dim_{{\mathbb Q}}\Gamma$. We assume that it is finite. 

Assume that the characteristics of a field $L$ and of the residue 
field $\kappa$ of a valuation $w$ of $L$ are equal. Then $w$ is called 
{\sl discrete}, if $L$ is algebraic over the subfield generated by 
a lift of a transcendence base of $\kappa$ and by a lift of a basis 
of the valuation group. In particular, $v$ is discrete if and only 
if the transcendence degree of $F$ over $F'$ is equal to $r$. 

The basic example of ${\mathcal O}_v$ is, after choosing an arbitrary 
algebraically closed $F'\subseteq F$, over which $F$ is of transcendence 
degree $r$, a transcendence base $x_1,\dots,x_r$ of $F$ over $F'$, and 
embeddings $F_j\hookrightarrow\lim\limits_{_N\longrightarrow\phantom{_N}}
F_{j-1}((x_j^{1/N}))$ over $F_{j-1}(x_j)$, the preimage in $F$ of the 
ring $\widehat{{\mathcal O}}_r=F'\oplus\widehat{{\mathfrak p}}_1=
\lim\limits_{_N\longrightarrow\phantom{_N}}F'[[x_1^{1/N}]]\oplus
\widehat{{\mathfrak p}}_2=\widehat{{\mathcal O}}_{j-1}\oplus
\widehat{{\mathfrak p}}_j$. Here $F_0:=F'$ and $F_j=\overline{F_{j-1}(x_j)}$, 
$\widehat{{\mathfrak p}}_{r+1}=0$ and $\widehat{{\mathfrak p}}_j=
\lim\limits_{_N\longrightarrow\phantom{_N}}x_j^{1/N}F'((x_1^{1/N}))
\dots((x_{j-1}^{1/N}))[[x_j^{1/N}]]\oplus\widehat{{\mathfrak p}}_{j+1}$ 
are prime ideals for all $1\le j\le r+1$. In this case 
$v(x_1^{m_1})<\dots<v(x_r^{m_r})$ for all $m_1,\dots,m_r>0$.

\label{obratim-endomor} If $r<\infty$ and $\sigma({\mathcal O}_v)\subseteq
{\mathcal O}_v$ for some $\sigma\in G$ then $\sigma\in G_v$, since 
$\sigma$ induces surjective endomorphism of $\Gamma$, i.e. an automorphism. 

It is well-known, \cite{zar-sam}, or Exercise 32, Chapter 5, \cite{am}, that 
the complement ${\mathcal O}\smallsetminus{\mathfrak p}$ to a prime ideal 
${\mathfrak p}$ in a valuation ring ${\mathcal O}$ projects onto the set 
of all non-negative elements in some isolated subgroup of the valuation 
group $\Gamma$, and ${\bf Spec}{\mathcal O}\stackrel{\sim}{\longrightarrow}
\{\mbox{isolated subgroups in $\Gamma$}\},$ ${\mathfrak p}\mapsto\langle 
v({\mathcal O}\smallsetminus{\mathfrak p})\rangle$, so there are exactly 
$r+1$ prime ideals in ${\mathcal O}_v$. 

{\sc Remarks.} 1. If ${\mathfrak p}\neq 0$ is a non-maximal prime ideal 
of finite codimension in ${\mathcal O}_v$ and ${\mathcal O}_{v'}:=
({\mathcal O}_v)_{{\mathfrak p}}$ then $G_v\subseteq G_{v'}$ (since any 
element $\sigma\in G_v$ preserves ${\mathfrak p}$, thus also 
${\mathcal O}_v\smallsetminus{\mathfrak p}$, i.e. induces 
an automorphism of ${\mathcal O}_{v'}$).
 
2. The inclusion 
$G_v\subset G_{\{F,{\mathcal O}_v[x_1^{-1}]\}|k}$ is proper for $r>1$, 
i.e. $G_v$ is not maximal. 

\vspace{4mm}

Let ${\mathcal P}^r_L$ be the set of discrete valuation rings of rank 
$r$ in $L$, containing $k$, admitting also the following description. 
Let ${\mathcal C}^r_X$ be the set of chains of irreducible 
normal subvarieties up to codimension $r$ on an irreducible proper
normal variety $X$ over $k$. Any proper surjection with irreducible 
fibres, e.g. a birational morphism, $X'\stackrel
{\pi}{\longrightarrow}X$ induces an embedding ${\mathcal C}^r_X
\hookrightarrow{\mathcal C}^r_{X'}$, $(Z^1\supset\dots\supset Z^r)
\mapsto(W^1\supset\dots\supset W^r)$, where $W^0:=X'$ and $W^j:=
(\pi|_{W^{j-1}})^{-1}_{{\rm prop}}(Z^j)$ for $1\le j\le r$ 
(and $\pi|_{W_1}:W_1\stackrel{\sim}{\longrightarrow}Z_1$). 
If $L$ is of finite type over $k$ then ${\mathcal P}^r_L\cong
\lim\limits_{_X\longrightarrow\phantom{_X}}{\mathcal C}^r_X$, 
where $X$ runs over the models of $L|k$, and ${\mathcal P}^r_F=
\left(\lim\limits_{\phantom{_L}\longleftarrow_L}(\coprod_{j=0}^r
{\mathcal P}^j_L)\right)\smallsetminus\coprod_{j=0}^{r-1}{\mathcal P}^j_F$. 
In particular, ${\mathcal P}^1_L=C(k)$, if $n=1$, where 
$C$ is a smooth proper model of the field $L$ over $k$. 

Any proper surjection $X'\stackrel{\pi}{\longrightarrow}X$ induces 
embeddings ${\mathbb Z}[{\mathcal C}^r_X]\hookrightarrow{\mathbb Z}
[{\mathcal C}^r_{X'}]$ and ${\mathbb Z}[{\mathcal C}^r_X]^{\circ}
\hookrightarrow{\mathbb Z}[{\mathcal C}^r_{X'}]^{\circ}$, where 
${\mathbb Z}[{\mathcal C}^r_X]^{\circ}:=\bigcap_{j=0}^{r-1}\ker
\left({\mathbb Z}[{\mathcal C}^r_X]\longrightarrow{\mathbb Z}
[{\mathcal C}^{r,j}_X]\right)$, ${\mathcal C}^{r,j}_X$ denotes 
the set of chains with no component of codimension $j$ and 
${\mathcal C}^r_X\longrightarrow{\mathcal C}^{r,j}_X$ is 
the omitting of such component. 

In particular, one can define a smooth $G$-module 
{\fontencoding{OT2}\selectfont C}$^r=\lim\limits
_{_L\longrightarrow\phantom{_L}}{\mathbb Z}[{\mathcal P}^r_L]^{\circ}$, 
where $L\subset F$ runs over the set of subfields of finite type over $k$. 
Then one can define a morphism $gr^W_{2q}H^q_{{\rm dR}/k}(F)
\stackrel{{\rm Res}}{\longrightarrow}${\fontencoding{OT2}\selectfont C}$^q
\otimes k$ by $\alpha\frac{dt_1}{t_1}\wedge\dots\wedge\frac{dt_q}{t_q}
\mapsto\sum_{\sigma\in{\mathfrak S}_q}{\rm sgn}(\sigma)\alpha|_{D_{1\dots q}}
\cdot(D_{\sigma(1)}\supset D_{\sigma(1)\sigma(2)}\supset\dots\supset
\bigcap_{j=1}^qD_j=:D_{1\dots q})$, where $D_i$ is given locally by 
$t_i=0$ and $\alpha$ is regular in a neighbourhood of $t_1=\dots=t_q=0$. 

\vspace{4mm}

The following fact is a well-known property of the Newton polygon. 
\begin{lemma} \label{newton} Let ${\mathcal O}$ be a valuation ring in 
an algebraically closed field $L$ and $v$ be the valuation. Suppose 
that the coefficients of a polynomial $P(X)=\sum_sa_sX^s\in L[X]$ 
satisfy $a_ia_j\neq 0$ for some $i<j$ and $(j-i)\cdot v(a_s)\ge
(j-s)\cdot v(a_i)-(i-s)\cdot v(a_j)$ for any $s$. 

Then there is $z\in L$ such that $P(z)=0$ 
and $(j-i)\cdot v(z)=v(a_i/a_j)$. \end{lemma}
\begin{proof} Replacing $a_s$ by $a_s\cdot a_i^{-1}\cdot t^{s-i}$, 
where $t^{j-i}=a_i/a_j$, and $z$ by $z/t$, we may assume that 
$a_s\in{\mathcal O}$ for any $s$ and $a_i,a_j\in{\mathcal O}^{\times}$, 
and look for $z\in{\mathcal O}^{\times}$ such that $P(z)=0$. 

If there is no such $z$ then $P(X)=cQ(X)R(X)$, where $c\in L^{\times}$, 
$Q(X)=\prod_{\alpha\in S}(X-\alpha)$, $R(X)=\prod_{\beta\in T}
(1-\beta X)$ and $S,T\subset{\mathfrak m}:={\mathcal O}\smallsetminus
{\mathcal O}^{\times}$. Thus, modulo ${\mathfrak m}$, $Q(X)R(X)$ 
is a power of $x$, contradicting the existence of two non-zero 
(modulo ${\mathfrak m}$) coefficients of $P(X)$. \end{proof} 

\begin{lemma} \label{transitiv} If $0\le r<n+1\le\infty$ then 
the $G$-action on the set of pairs $(v,\Lambda)$, where $v:F^{\times}
\longrightarrow\!\!\!\!\to\Gamma\cong{\mathbb Q}^r$ is an element of 
${\mathcal P}^r_F$ and $\Lambda\cong{\mathbb Z}^r$ is a lattice in 
$\Gamma$, is transitive. The stabilizer of $(v,\Lambda)$ acts 
transitively on the set of maximal subfields $\widetilde{F}$ in $F|k$ 
such that $v(\widetilde{F}^{\times})=\Lambda$. The residue field of 
$\widetilde{F}$ coincides with $\kappa(v)$ (in particular, it is 
algebraically closed). \end{lemma} 
\begin{proof} Let $0=\Gamma_0\subset\Gamma_1\subset\dots\subset
\Gamma_{r-1}\subset\Gamma_r=\Gamma$ be the isolated subgroups in 
$\Gamma$. Choose $x_1,x_2,\dots\in{\mathcal O}_v$ such that $v(x_j)
\in\Gamma_j\smallsetminus\Gamma_{j-1}$ for $1\le j\le r$ and 
$x_{r+1},x_{r+2},\dots$ modulo $\mathfrak{m}_v$ form a transcendence 
base of ${\mathcal O}_v/\mathfrak{m}$ over $k$. Clearly, 
$x_1,x_2,\dots$ are algebraically independent over $k$. 
Set $k'=\overline{k(x_{r+1},x_{r+2},\dots)}\subseteq F$ and 
$\widehat{F}_r:=\lim\limits_{_N\longrightarrow\phantom{_N}}
k'((X_1^{1/N}))\dots((X_r^{1/N}))$.

Define the embedding $k(x_1,x_2,\dots)\stackrel{\varphi}{\hookrightarrow}
k'((X_1))\dots((X_r))$ by $x_j\mapsto X_j$. It respects the valuation. 
Consider the set $S$ of embeddings into $\widehat{F}_r$ of subfields in 
$F$ containing $k(x_1,\dots,x_n)$, extending $\varphi$ and respecting 
the valuation. The set $S$ is partially ordered. Clearly, $S$ contains 
maximal elements. Let $L'\stackrel{\xi}{\hookrightarrow}\widehat{F}_r$ 
be a maximal element of $S$. If $\sum_sa_sy^s=0$ is a minimal polynomial 
of $y\in F$ over $L'$ then there exist $i<j$ such that $v(y)=
\frac{v(a_i/a_j)}{j-i}$. By Lemma \ref{newton}, there exists 
$z\in\widehat{F}_r$ such that $\sum_s\xi(a_s)z^s=0$ and $\hat{v}(z)
=\frac{\hat{v}(\xi(a_i/a_j))}{j-i}=\frac{v(a_i/a_j)}{j-i}$. 
Therefore, $\xi$ extends to $L'(z)$ by $y\mapsto z$, and thus, $L'=F$. 

This shows that $\varphi$ extends to an embedding 
$F\hookrightarrow\widehat{F}_r$, respecting the valuation. 

Thus, any element of ${\mathcal P}^r_F$ is determined by a choice 
of $x_1,x_2,\dots\in F$ algebraically independent over $k$ and by 
an embedding of $F$ into $\widehat{F}_r$ over $k(x_1,x_2,\dots)$. 
Clearly, the $G$-action is transitive on these data, which means its 
transitivity on ${\mathcal P}^r_F$. 

The completion of $\widetilde{F}$ is isomorphic to $k'((X_1))\dots((X_r))$, 
so $\widetilde{F}$ is isomorphic to the algebraic closure of 
$k(x_1,x_2,\dots)$ in $k'((X_1))\dots((X_r))$. As $G_v$ acts transitively 
on the lattices in $\Gamma$, the isomorphism class of the extension 
$\widetilde{F}|k$ is determined uniquely. \end{proof} 

\vspace{4mm}

The $G_v$-action on $\kappa(v)$ induces a homomorphism 
$G_v\longrightarrow\hspace{-3mm}\to G_{\kappa(v)|k}$. Let 
\begin{multline*} G_v^{\dagger}:=\{\sigma\in G_v~|~\mbox{$\sigma x-x
\in{\mathfrak m}_v$ for any $x\in{\mathcal O}_v$}\}\\ 
=\{\sigma\in G~|~\mbox{$\frac{\sigma x}{x}\in 1+{\mathfrak m}_v$ 
for any $x\in{\mathcal O}_v^{\times}$}\}\end{multline*} 
be its kernel, the `inertia' subgroup. 

A continuous\footnote{One can apply the well-known fact (cf., e.g., \cite{W}) 
that the ${\mathbb Z}$-subalgebra in $k$ generated by the coefficients of 
any element of $k((t))$, which is algebraic over $k(t)$, has a finite number 
of generators. {\it Proof.} Let $F(f,t)=0$ be a minimal polynomial of some 
formal series $f\in k[[t]]$. Then the first partial derivative 
of $F$ with respect to first variable does not vanish at $(f,t)$: 
$F_I(f,t)\neq 0$, so $F_I(f,t)\in t^sk[[t]]\smallsetminus
t^{s+1}k[[t]]$ for some $s\ge 0$. 

Denote by $f_n$ the only polynomial of degree $<n$ congruent to 
$f$ modulo $t^n$. For a polynomial $\Phi$ in $t$ and an integer 
$m$ denote by $\Phi_{[m]}$ the degree-$m$ coefficient of $\Phi$. 
Clearly, $F_I(f_n,t)_{[m]}$ is independent of $n$ for $n>m$. 

We are going to show that the ${\mathbb Z}$-subalgebra in $k$ generated 
by the coefficients of $f$ is generated, in fact, by the coefficients of 
$F$, by the coefficients of $f_{s+1}$ and by the inverse of a polynomial 
in coefficients of $F$ and in coefficients of $f_{s+1}$. 

This is done by induction on degree $n>s$: by definition, 
$F(f_n,t)\in t^nk[[t]]$ and we have to find (assuming that it exists!) 
an element $a_n\in k$ such that $F(f_n+a_nt^n)\in t^{n+1}k[[t]]$. 
One has $F(f_n+a_nt^n,t)\equiv F(f_n,t)+F_I(f_n,t)a_nt^n\pmod{t^{2n}}$, 
so the condition is $F(f_n,t)_{[n+s]}+F_I(f_n,t)_{[s]}a_n=0$. 
This is a linear equation with polynomial coefficients in 
coefficients of $F$ and in coefficients of $f_{s+1}$. As $a_n$ is a 
(unique!) solution of this linear equation, $F_I(f_n,t)_{[s]}=:D$ 
is non-zero. Then the coefficients $a_n$ are polynomials over ${\mathbb Z}$ 
in coefficients of $F$, coefficients of $f_{s+1}$ and in $D^{-1}$. \qed} 
section of $G_{\kappa(v)|k}\hookrightarrow G_v$ is determined by a 
subfield $F'\subseteq{\mathcal O}_v$, identified with $\kappa(v)$ by 
the reduction modulo $\mathfrak{m}_v$, and by an embedding of $F$ 
into the field of iterated Puiseux series 
$\lim\limits_{_N\longrightarrow\phantom{_N}}F'((t_1^{1/N}))
\dots((t_r^{1/N}))$ over $F'$ compatible with valuations 
(i.e. by a choice of a section $\Gamma\hookrightarrow F^{\times}$ 
of the valuation $v$). 

As the centralizer of $G^{\dagger}_v$ is trivial, such sections form 
a $G^{\dagger}_v$-torsor. One has a decomposition $G_v=G_v^{\dagger}
(G_v\cap G_{\{F,F'\}|k})$, where $G^{\dagger}_v\cap G_{\{F,F'\}|k}=
G_v\cap G_{F|F'}$. 

Restriction to $F'$ and reduction modulo ${\mathfrak m}_v$, 
respectively, determine canonical isomorphisms $G_{F'|k}
\stackrel{\sim}{\longleftarrow}(G_v\cap G_{\{F,F'\}|k})/
(G_v\cap G_{F|F'})\stackrel{\sim}{\longrightarrow}G_{\kappa(v)|k}$. 

\vspace{4mm}

Let $B_+$ be the group of linear transformations of $\Gamma$ 
(isomorphic to the group of upper-triangular rational $r\times r$ 
matrices with positive diagonal entries), respecting the order, i.e. 
the flag of isolated subgroups $0=\Gamma_0\subset\Gamma_1\subset\dots
\subset\Gamma_r=\Gamma$ and the orientation of each $\Gamma_j/
\Gamma_{j-1}\cong{\mathbb Q}$, $1\le j\le r$. Then 
${\mathcal P}^r_L$ consists of (some) valuations $(L^{\times}
/k^{\times})\otimes{\mathbb Q}\longrightarrow\!\!\!\!\!\rightarrow
{\mathbb Q}^r$ modulo the $B_+$-action. The $G_v$-action on 
$F^{\times}$ induces a linear action on $\Gamma$, preserving the 
order, i.e. a homomorphism $G_v\longrightarrow B_+$, surjective 
and split, if $v$ is discrete. 

Let $G_v^{\circ}:=\{\sigma\in G_v~|~\mbox{$\sigma x/x\in
{\mathcal O}_v^{\times}$ for any $x\in F^{\times}$}\}$ be its kernel.  

There is a homomorphism $G_v^{\dagger}\cap G_v^{\circ}\longrightarrow
{\rm Hom}(\Gamma,\kappa(v)^{\times})$,\footnote{Clearly, ${\rm Hom}
(\Gamma,\kappa(v)^{\times})\cong(\!\!\lim\limits_{\phantom{_N}
\longleftarrow_N}M_N)^r$, where $M_{NN'}:=\kappa(v)^{\times}
\longrightarrow\hspace{-2mm}\to M_{N'}:=\kappa(v)^{\times}$ 
is the raising to the $N$-th power.} 
$\sigma\longmapsto(v(x)\mapsto\sigma x/x\bmod\mathfrak{m}_v)$, 
which is surjective and split, if $v$ is discrete. 
\begin{lemma} \label{disre-plotno} Its kernel $G_v^1:=\{\sigma\in
G_v~|~\mbox{$\sigma x/x\in 1+\mathfrak{m}_v$ for any $x\in F^{\times}
$}\}\subset G_v^{\dagger}$ is a discrete subgroup if $n<\infty$. The 
functors $H^0(G,-)$ and $H^0(G^1_v,-)$ coincide on the category of smooth 
$G$-sets, if $n=\infty$ and $v$ is discrete of finite non-zero rank. 
\end{lemma}
\begin{proof} The $G_v$-action on $F$ extends to a continuous 
$G_v$-action on the completion $F_v$ of $F$, and the continuous 
$G^1_v$-action on $F_v$ is determined by the action on some $r$-tuple 
of elements, representing the isolated subgroups of $\Gamma$ and 
by the action on some maximal subfield with trivial restriction of $v$. 
Clearly, the latter action is determined by its restriction to any
transcendence base over $k$. The second assertion follows from the 
fact that any open subgroup containing $G^1_v$ coincides with $G$, 
which follows from Proposition \ref{morfizm-mono}. \end{proof}

\vspace{4mm}

{\sc Remark.} There is a natural inclusion of sets $G_v^{\circ}/G_v^1
\hookrightarrow{\rm Hom}(F^{\times}/(k^{\times}+\mathfrak{m}_v),
\kappa(v)^{\times})$, $[\sigma]\mapsto(\sigma x/x\bmod\mathfrak{m}_v)$. 

Let $L$ be the function field of a $d$-dimensional variety over $k$, 
$I\subseteq\{1,\dots,r\}$ be a subset and $v\in{\mathcal P}^r_F$, 
$p\in{\mathcal P}^{|I|}_L$. Let $O_{p,v,I}$ be the set of all 
embeddings $\sigma:L\stackrel{/k}{\hookrightarrow}F$ such that 
$\sigma^{-1}({\mathcal O}_v)={\mathcal O}_p$ and 
$\sigma(L^{\times})\cap\Gamma_i\neq\sigma(L^{\times})\cap\Gamma_{i-1}$ 
if and only if $1\le i\le r$ and $i\in I$. 
\begin{proposition} \label{specializ} If $m:=\max(0,r+d-n)\le|I|\le 
M:=\min(d,r)$ then $O_{p,v,I}$ is a non-empty $G_v$-orbit. The set 
$\{L\stackrel{/k}{\hookrightarrow}F\}$ 
is a disjoint union of $O_{p,v,I}$. In particular, 
${\mathbb Q}[\{L\stackrel{/k}{\hookrightarrow}F\}]_{G^{\dagger}_v}=
\bigoplus_{s=m}^M\bigoplus_{p\in{\mathcal P}^s_L}
{\mathbb Q}[\{\kappa(p)\stackrel{/k}{\hookrightarrow}\kappa(v)\}]
^{\binom{r}{s}}$ and $G_v\backslash\{L\stackrel{/k}{\hookrightarrow}
F\}=\coprod_{s=m}^M({\mathcal P}^s_L)^{\coprod\binom{r}{s}}$. 
\end{proposition}
\begin{proof} Any element of $\sigma(L^{\times})$ is either contained 
in ${\mathcal O}_v$, or its inverse is contained in ${\mathcal O}_v$, 
i.e., ${\mathcal O}_p:=\sigma^{-1}({\mathcal O}_v)$ is a valuation 
ring.\footnote{Note, that the henselian property does not suffice. 
For instance, if $k(X)=k(x_1,\dots,x_d)$ is embedded into the field 
of fractions of ${\mathcal O}:=k[[T_1,\dots,T_N]]$, $x_j\mapsto 
T^{\alpha_j}$, $|\alpha_j|=0$ then $k(X)\cap{\mathcal O}=k$, since 
it consists of homogeneous functions in $T_1,\dots,T_N$ of degree 0.} 
Clearly, ${\mathcal O}_p$ is a discrete valuation ring. 

As $\sigma:L^{\times}/{\mathcal O}_p^{\times}\hookrightarrow F^{\times}
/{\mathcal O}_v^{\times}$, one has ${\rm rk}(L^{\times}/{\mathcal O}_p
^{\times})\le M$. Let us check that ${\rm rk}(L^{\times}/{\mathcal O}_p
^{\times})\ge m$, and in particular, the inclusion $(L^{\times}/
{\mathcal O}_p^{\times})\otimes{\mathbb Q}\stackrel{\sigma}{\hookrightarrow}
\overline{\sigma(L^{\times})}/\overline{\sigma(L^{\times})}\cap{\mathcal O}
_v^{\times}\subseteq\Gamma$ is bijective. Let $x\in F^{\times}$ be 
algebraic over $\sigma(L)$ and $\sum_{j=0}^ma_jx^j=0$ 
be a minimal polynomial of $x$ over $\sigma({\mathcal O}_p)$. Then there 
exist $0\le i<j\le m$ such that $a_ia_j\neq 0$ and $v(a_ix^i)=v(a_jx^j)$, 
so $v(x)=\frac{v(a_i/a_j)}{j-i}$, i.e., the element $v(x)\in 
F^{\times}/{\mathcal O}_v^{\times}$ is in the image of $(L^{\times}
/{\mathcal O}_p^{\times})\otimes\frac{1}{m!}{\mathbb Z}$.

We may assume that $a_j=1$ for some $j$. Thus, if $x\in
{\mathcal O}_v$ then its image in ${\mathcal O}_v/\mathfrak{m}_v$ 
is algebraic over $\sigma({\mathcal O}_p)/\sigma(\mathfrak{m}_p)$, 
i.e. ${\mathcal O}_v/\mathfrak{m}_v\cong\overline{{\mathcal O}_p
/\mathfrak{m}_p}$, and therefore, $p\in{\mathcal P}^r_L$, if $d=n$. 

To check that $O_{p,v,I}$ is a $G_v$-orbit, we have to show that for 
any pair of its elements $\sigma,\tau$ there exists $\xi\in G_v$, 
extending $\tau\sigma^{-1}:\sigma({\mathcal O}_p)\stackrel{\sim}
{\longrightarrow}\tau({\mathcal O}_p)$. The embeddings $\sigma,\tau:
L\stackrel{/k}{\hookrightarrow}F$ are such that $\sigma({\mathcal O}_p),
\tau({\mathcal O}_p)\subset{\mathcal O}_v$ and $\sigma(L)^{\times}$ and 
$\tau(L)^{\times}$ span some subspaces in $\Gamma$ in the same position 
with respect to the flag of isolated subgroups $\Gamma_1\subset\Gamma_2
\subset\dots\subset\Gamma_r=\Gamma$. 

The set $$S:=\left\{A\stackrel{\varphi}{\hookrightarrow}{\mathcal O}_v
~\bigg\vert~\begin{array}{l}\mbox{$v(\varphi(x))\ge v(\varphi(y))$,
if $v(x)\ge v(y)$ for $x,y\in A$,}\\
\sigma({\mathcal O}_p)\subseteq A\subseteq{\mathcal O}_v,\quad
\varphi|_{\sigma({\mathcal O}_p)}=\tau\sigma^{-1}\end{array}
\right\}$$ is non-empty, since it contains 
$(\sigma({\mathcal O}_p)\stackrel{\tau\sigma^{-1}}
{\hookrightarrow}{\mathcal O}_v)$. According to Zorn's lemma, 
there are maximal elements in $S$, for instance, 
$B\stackrel{\psi}{\hookrightarrow}{\mathcal O}_v$. Then $B$ is 
integrally closed in ${\mathcal O}_v$, since for any 
$x\in{\mathcal O}_v$ integral over $B$ an embedding $\psi$ 
extends to $B[x]\subseteq{\mathcal O}_v$, even respecting the 
order,\footnote{i.e., if $\sum_{s=0}^ma_sx^s=0$, $a_s\in B$ 
and $v(x)=\frac{v(a_i/a_j)}{j-i}$ for some $i<j$ then one can choose 
$y\in{\mathcal O}_v$ such that $\sum_{s=0}^m\psi(a_s)y^s=0$ and 
$v(y)=\frac{v(\psi(a_i/a_j))}{j-i}$.} since ${\mathcal O}_v$ is 
integrally closed in $F$. Thus, any element of ${\mathcal O}_v$ can 
be presented as $a/b$, where $a,b\in B$ and $v(a)\ge v(b)$. As $\psi$ 
respects the order, it maps ${\mathcal O}_v$ into ${\mathcal O}_v$, i.e., 
${\mathcal O}_v=B$ and $\psi\in G_v$, cf. p.\pageref{obratim-endomor}. 

The set $O_{p,v,I}$ is non-empty, since, given an embedding 
$L\stackrel{/k}{\hookrightarrow}F$, the valuation $p$ extends to 
an element of ${\mathcal P}^r_F$, and the group $G$ permutes the 
elements of ${\mathcal P}^r_F$ (Lemma \ref{transitiv}). \end{proof}

\subsection{Valuations and maximal subgroups}
{\sc Example.} If $n=1$ then to any valuation $v$ the decomposition 
$\{k(C)\stackrel{/k}{\hookrightarrow}F\}=C(F)\smallsetminus C(k)=
\coprod_{C(k)}({\mathfrak m}_v\smallsetminus\{0\})$ is associated. 

\begin{proposition} \label{peresech-orbit} If $n=1$ then $O_{p,v}
\cap O_{q,v'}$ is non-empty for any pair of distinct $v,v'\in
{\mathcal P}^1_F$ and any pair of distinct points $p,q\in C(k)$ 
on a smooth proper curve $C$ over $k$, i.e. there is an embedding 
$\sigma:k(C)\stackrel{/k}{\hookrightarrow}F$ such that $\sigma
({\mathcal O}_p)\subset{\mathcal O}_v$ and $\sigma({\mathcal O}_q)
\subset{\mathcal O}_{v'}$. (Here $O_{p,v}:=O_{p,v,\{1\}}$.) \end{proposition}
\begin{proof} We consider $v$ and $v'$ as a compatible system of 
points on smooth proper curves over $k$ with the function fields 
embedded into $F$: if $C_{\beta}\longrightarrow\!\!\!\!\to 
C_{\alpha}$ then $v_{\beta}\mapsto v_{\alpha}$ and 
$v'_{\beta}\mapsto v'_{\alpha}$. One needs to check that there exist 
$\beta$ and a map $C_{\beta}\longrightarrow\!\!\!\!\to C$ such that 
$v_{\beta}\mapsto p$ and $v'_{\beta}\mapsto q$. 

Choose a non-constant function $x\in{\mathcal O}(C\smallsetminus\{p\})$, 
i.e. a surjective morphism $C\longrightarrow{\mathbb P}^1$ sending 
to $\infty$ only $p$. Let $C'\longrightarrow C$ be a cover such that 
the composition $C'\longrightarrow C\longrightarrow{\mathbb P}^1$ is 
a Galois cover with the group $A$. For some $\alpha$ such that 
$v_{\alpha}\neq v'_{\alpha}$ choose a surjection $C_{\alpha}
\longrightarrow{\mathbb P}^1$ such that $v_{\alpha}\mapsto\infty$ 
and $v'_{\alpha}\mapsto x(q)$. Consider the normalization $D$ of 
an irreducible component of $C'\times_{{\mathbb P}^1}C_{\alpha}$.
The surjection $D\longrightarrow C_{\alpha}$ is isomorphic to the 
surjection $C_{\beta}\longrightarrow C_{\alpha}$ for some $\beta$.
Let $\pi:C_{\beta}\longrightarrow C'$ be the projection. As $A$ 
acts transitively on the fibres of the composition $C'\longrightarrow
C\longrightarrow{\mathbb P}^1$, there is an element $\gamma\in A$ 
such that $\gamma\pi(v'_{\beta})$ belongs to the preimage of $q$. 
Then the composition $C_{\beta}\stackrel{\gamma\pi}{\longrightarrow}C'
\longrightarrow C$ maps $v_{\beta}$ to $p$, and $v'_{\beta}$ to $q$. \end{proof}

\begin{proposition} For any pair of distinct $v,v'\in{\mathcal P}^1_F$ 
the subgroup $H$, generated by $G_v$ and $G_{v'}$, acts transitively 
on $\{L\stackrel{/k}{\hookrightarrow}F\}$, i.e. $H$ is dense in $G$. 
\end{proposition} 
\begin{proof} Let $L$ be an extension of $k$ and $w,w':L^{\times}/
k^{\times}\stackrel{\not\equiv 0}{\longrightarrow}{\mathbb Q}$ be 
a pair of discrete valuations. If ${\mathcal O}^{\times}_w\subseteq
{\mathcal O}^{\times}_{w'}$ in $L$ then ${\mathcal O}^{\times}_{w'}
/{\mathcal O}^{\times}_w$ is the kernel of ${\mathbb Q}\supseteq 
L^{\times}/{\mathcal O}^{\times}_w\longrightarrow\hspace{-3mm}\to 
L^{\times}/{\mathcal O}^{\times}_{w'}\subseteq{\mathbb Q}$, which is 
evidently injective, so $w=w'$. Let $w\neq w'$. Then for any $x\in
{\mathcal O}^{\times}_w$ one has either $x\in{\mathcal O}_{w'}$, 
or $x^{-1}\in{\mathcal O}_{w'}$, that is there is some $t\in
{\mathcal O}^{\times}_w\cap{\mathcal O}_{w'}$ such that 
$t\not\in{\mathcal O}_{w'}^{\times}$. Fix such $t$ and $t^{1/N}$ 
for all integers $N\ge 1$. Let $x_1,x_2,\dots\in{\mathcal O}^{\times}_w$ 
be a lifting of a transcendence base of $\kappa(w)$ over $k(t)$. 
Set $k_0:=k$ and define inductively a strictly increasing 
sequence of algebraically closed subfields in 
$\{0\}\cup({\mathcal O}^{\times}_w\cap{\mathcal O}^{\times}_{w'})
\subset L$ as follows. For any $i\ge 1$ there exist $P\in k_{i-1}[T]$, 
an integer $N\ge 1$ and $M\in{\mathbb Q}^{\times}_+$ such that $y_i:=
t^{-M}(x_i^{w'(t)}-t^{w'(x_i)}P(t^{1/N}))\in{\mathcal O}^{\times}_{w'}$ 
and $y_i\not\in k_{i-1}+{\mathfrak m}_{w'}$. Then $y_1,y_2,\dots$ is 
another lifting of a transcendence base of $\kappa(w)$ over $k(t)$ 
in ${\mathcal O}^{\times}_w$. Set $k_i:=\overline{k_{i-1}(y_i)}$ and 
$k'=\bigcup_{i\ge 1}k_i$. Then $\kappa(w)$ is algebraic over the 
reduction of $k'$ modulo ${\mathfrak m}_w$ and $\kappa(w')$ is 
algebraic over the reduction of $k'$ modulo ${\mathfrak m}_{w'}$. 
This shows that we are reduced to the case of $n=1$. 

By Proposition \ref{peresech-orbit}, for any $\xi:L
\stackrel{/k}{\hookrightarrow}F$ the map $G_v\times G_{v'}
\longrightarrow\{L\stackrel{/k}{\hookrightarrow}F\}$, given 
by $(\sigma,\tau)\mapsto\sigma\tau\xi$, is surjective. \end{proof} 

\vspace{4mm}

If $n=r=1$ define $\varphi:G_v^1\longrightarrow\Gamma\cup\{+\infty\}$ 
by $\sigma\mapsto v(\sigma x/x-1)$ for any $x\in{\mathfrak m}
\smallsetminus\{0\}$, or $x\in F\smallsetminus{\mathcal O}_v$. 
\begin{lemma} $\varphi$ is independent of $x$ and determines 
a bounded non-archimedian bi-invariant distance on $G_v^1$. The 
logarithmic distance transforms the adjoint $G_v$-action on $G_v^1$ 
to the natural $G_v$-action on $\Gamma\cong{\mathbb Q}$. \end{lemma}
\begin{proof} Independence of $x$: if $\sigma\neq 1$, $y=
\sum_{j\ge 1}a_jx^{j/N}\neq 0$ for some $a_j\in k$ and $\sigma x/x=
1+ax^{\alpha}$, where $a\in{\mathcal O}_v^{\times}$ and $\alpha>0$, 
then $\sigma y-y=\sum_{j\ge 1}a_jx^{j/N}((1+ax^{\alpha})^{j/N}-1)$, 
i.e. \begin{equation}\label{odnoiz}\sigma y/y-1=
\frac{v(y)}{v(x)}ax^{\alpha}+o(x^{\alpha}).\end{equation}
Thus, $v(\sigma y/y-1)=v(x^{\alpha})=v(\sigma x/x-1)$.

From the equality $\sigma\tau x/x-1=(\sigma(\tau x)/\tau x-1)\tau x/x+
\tau x/x-1$ we get that $\varphi(\sigma\tau)\ge\min(\varphi(\sigma),
\varphi(\tau))$, $\varphi(\xi\tau\xi^{-1})=\xi\varphi(\tau)$ for any 
$\xi\in G_v$ (in particular, $\varphi(\xi\tau\xi^{-1})=\varphi(\tau)$ 
for any $\xi\in G_v^{\circ}$) and $\varphi(\sigma)=\varphi(\sigma^{-1})$:
$v(\sigma^{-1}x/x-1)=v(x/\sigma x-1)=\varphi(\sigma)$. \end{proof} 

\begin{lemma} For any integer $N\ge 1$ the self-map of $G_v^1$, 
$\sigma\mapsto\sigma^N$, is injective. \end{lemma}
\begin{proof} This follows from the identity $(f+h)^{\circ N}=
f^{\circ N}+Nh+o(h)$ for any formal Puiseux series $f=x+o(x)$ 
and $h=o(x)$ which can be checked by induction on $N$: for $N=1$ 
this is evident; $(f+h)^{\circ (N+1)}=(f+h)(f^{\circ N}+Nh+o(h))
=f(f^{\circ N}+Nh+o(h))+h+o(h)=f^{\circ (N+1)}+f'(f^{\circ N})
\cdot Nh+o(h)+h+o(h)=f^{\circ (N+1)}+(N+1)\cdot h+o(h)$. \end{proof} 

\vspace{4mm}

Let $G^1_v(\beta):=\{\sigma\in G^1_v~|~\varphi(\sigma)\ge\beta\}$,
where $\beta\in\Gamma\otimes{\mathbb R}$. This is a normal subgroup 
in $G^{\circ}_v$. Then $G^1_v=G^1_v(0)=G^1_v(0)^+$, where 
$G^1_v(\beta)^+:=\bigcup_{\gamma>\beta}G^1_v(\gamma)=
\{\sigma\in G^1_v~|~\varphi(\sigma)>\beta\}$. Clearly, $G^1_v(\beta)
\neq G^1_v(\gamma)$, if $\beta\neq\gamma$.
\begin{lemma} There is a canonical isomorphism $G^1_v(\beta)
/G^1_v(\beta)^+\stackrel{\sim}{\longrightarrow}{\rm Hom}
(\Gamma,{\mathfrak m}^{[\beta]})$, where 
$${\mathfrak m}^{[\beta]}=\{x\in{\mathfrak m}~|~v(x)\ge\beta\}
/\{x\in{\mathfrak m}~|~v(x)>\beta\}\cong\left\{\begin{array}{ll}
k, & \mbox{if $\beta\in\Gamma$ and $\beta>0$}\\
0, & \mbox{otherwise}\end{array}\right.$$ \end{lemma}
\begin{proof} It is clear from the formula (\ref{odnoiz}) that any 
element $\sigma\in G^1_v(\beta)$ induces a homomorphism $\Gamma
\longrightarrow{\mathfrak m}^{[\beta]}$, $v(y)\mapsto\sigma y/y-1$.

This gives a surjective homomorphism $G^1_v(\beta)/G^1_v(\beta)^+
\longrightarrow{\rm Hom}(\Gamma,{\mathfrak m}^{[\beta]})$:
$\sigma\xi\longmapsto(v(x)\mapsto\sigma\xi x/x-1=
\sigma(\xi x/x-1)\sigma x/x+(\sigma x/x-1)
\equiv(\xi x/x-1)+(\sigma x/x-1))$.

Let us check its injectivity: $\sigma x=x(1+ax^{\beta}+x^{\beta}h)$, 
$h\in{\mathfrak m}$. Let $\tau[x(1+ax^{\beta}+x^{\beta}h)]:=
x(1+ax^{\beta})$. Then $\tau\sigma x/x=1+ax^{\beta}$ and, if we set 
$y:=x(1+ax^{\beta})$ then $x(1+ax^{\beta}+x^{\beta}h)=y+o(y^{\beta+1})$, 
and thus, $\tau\in G^1_v(\beta)^+$. Finally, $G^1_v(\beta)/G^1_v(\beta)^+=
{\rm Hom}(\Gamma,{\mathfrak m}^{[\beta]})$. \end{proof}

\begin{lemma} \label{sur-fil} $G^1_v(\beta)$ is surjective over 
$G/G^{\circ}$ for any $\beta\in\Gamma\otimes{\mathbb R}$.\end{lemma}
\begin{proof} As $G_v$ is maximal and does not contain $G^{\circ}$, 
it suffices to show that the subgroup generated by $G^1_v(\beta)$ 
and $G^{\circ}$ contains $G_v$. A choice of a section $s:\Gamma
\longrightarrow F^{\times}$ of the projection $v:F^{\times}
\longrightarrow\Gamma$ determines an additive section ${\rm Hom}
(\Gamma,k^{\times})\longrightarrow G_v^{\circ}\cap G^{\circ}$
of the projection $G_v^{\circ}\longrightarrow{\rm Hom}(\Gamma,k^{\times})$.
A choice of $x\in{\mathfrak m}\cap s(\Gamma)$ gives a bijection of sets 
$G^1_v$ and ${\mathfrak m}$: ${\mathfrak m}\ni m:x\mapsto x(1+m)$.
The elements $(x\mapsto(1+x^{\alpha})^{1/\alpha}-1)\in G_v\cap G^{\circ}$,
$\alpha\in{\mathbb Q}^{\times}_+$, determine a set-theoretic section 
of the projection $G_v\longrightarrow{\mathbb Q}^{\times}_+$.

It remains to check that the subgroup generated by $G^1_v(\beta)$ 
and $G^1_v\cap G^{\circ}$ contains $G^1_v$. Clearly, 
$x(1+ax^{\alpha})^{-1/\alpha}\in G^1_v\cap G^{\circ}$.
If an element of $G^1_v$ is presented by $h=x(1+ax^{\alpha}+o(x^{\alpha}))
\in x(1+{\mathfrak m})\cap k[[x^{1/N}]]$ then 
$h(1+a\alpha h^{\alpha})^{-1/\alpha}=x(1+o(x^{\alpha}))
\in x(1+{\mathfrak m})\cap k[[x^{1/N}]]$. So the composition of 
$h$ with an appropriate element of $G^1_v\cap G^{\circ}$ 
will be in $G^1_v(\beta)$. \end{proof}

\section{Valuations and associated functors} \label{assoc-funktory}
The aim of this section is to associate a sheaf in the smooth topology 
to any smooth representation of $G$. This can be achieved in the following 
way. For each smooth $k$-variety $X$, a scheme point $p\in X$ and an 
embedding $k(X_p)\stackrel{/k}{\hookrightarrow}F$ of the function field 
of the connected component $X_p\ni p$ of $X$, define a collection 
$J_{X,p}$ of subfields in $F$ and then set ${\mathcal W}_{X,p}:=
W^{G_{F|k(X_p)}}\cap\left(\sum_{F''\in J_{X,p}}W^{G_{F|F''}}\right)$. 
This is what is called ``globalization'', cf. \S\ref{globalizaciq}. 

Obviously, for any $q\ge 0$ we want to get the sheaf of differential 
$q$-forms $\Omega^q_{{\mathcal O}|k}$ from the smooth representation 
$\Omega^q_{F|k}$ of $G$. For the representation 
$\bigotimes_F^q\Omega^1_{F|k}$ there are more options: 
one `homotopy invariant' $\bigotimes_{{\mathcal O}}^q
\Omega^1_{{\mathcal O}|k}$ (which is more reasonable), and another 
one with the Galois descent property. Here `homotopy invariant' means 
that for any projective bundle $X\longrightarrow Y$ with a proper 
base $Y$ the induced morphism of sections is an isomorphism. 

On the other hand, it is natural to ask, whether the constructed sheaves 
are functorial with respect to the regular embeddings. In some extremal 
situation this is discussed in \S\ref{specializaciq}. 

\subsection{The ``globalization'' functor} \label{globalizaciq} 
For any collection $J$ of subfields $F''\subset F$ the additive 
functor $$\Phi_J:W\mapsto\sum_{F''\in J}W^{G_{F|F''}}$$ on 
${\mathcal S}m_G$ preserves surjections if any element of $J$ is 
either of infinite transcendence degree over $k$, or is contained 
in an element of $J$ of arbitrarily big finite transcendence 
degree over $k$. In general, $\Phi_J$ preserves the injections. 

{\sc Remark.} If $J$ consists of all purely transcendental 
extensions of $k$ then $\Phi_J(\Omega^q_{F|k})=\Omega^q_{F|k}$ 
if $n>q$ and $\Phi_J(\Omega^q_{F|k,{\rm reg}})=0$ for any $q\ge 1$, 
so $\Phi_J$ is not exact in general, even if $n=\infty$. 

In particular, for a discrete valuation $v\in{\mathcal P}^r_F$ consider 
the functor $$(-)_v:{\mathcal S}m_G\longrightarrow{\mathcal S}m_{G_v},$$ 
$W\mapsto W_v:=\sum_{F''\subset{\mathcal O}_v}W^{G_{F|F''}}=
\sum_{\sigma\in G_v}W^{G_{F|\sigma(F')}}=
\sum_{\sigma\in G^{\dagger}_v}W^{G_{F|\sigma(F')}}\subseteq W$ for any 
algebraically closed subfield $F'\subset{\mathcal O}_v$ such that 
${\rm tr.deg}(F'|k)={\rm tr.deg}(F|k)-r$. 

Set $\Gamma_r(W):=\bigcap_{v\in{\mathcal P}^r_F}W_v$ and 
$\Gamma:=\Gamma_1$, so $\Gamma_r:{\mathcal S}m_G
\longrightarrow{\mathcal S}m_G$ are additive functors. 

{\sc Example.} ${\mathbb Q}[\{L\stackrel{/k}{\hookrightarrow}F\}]_v
={\mathbb Q}[\{L\stackrel{/k}{\hookrightarrow}{\mathcal O}_v\}]$ 
and $(F[\{L\stackrel{/k}{\hookrightarrow}F\}])_v={\mathcal O}_v
[\{L\stackrel{/k}{\hookrightarrow}{\mathcal O}_v\}]$ (and all these 
modules are zero if ${\rm tr.deg}(L|k)>n-r$). 

If $n=\infty$ then for any open subgroup $U\subseteq G$ the module 
${\mathbb Q}[G/U]_v$ (resp., $(F[G/U])_v$) is the direct sum of 
${\mathbb Q}\cdot[\sigma]$ (resp., of ${\mathcal O}_v\cdot[\sigma]$) over 
all $[\sigma]\in G/U$ such that $\sigma U\sigma^{-1}$ contains $G_{F|L}$ 
for some subfield $L\subset{\mathcal O}_v$. ({\it Proof.} There are no 
proper open subgroups of finite index in $G_{F|F''}$ for any algebraically 
closed subfield $F''$ over which $F$ is of infinite transcendence degree, 
cf. \cite{repr}. Thus, either $G_{F|F''}$ fixes the coset 
$[\sigma]$, or the $G_{F|F''}$-orbit of $[\sigma]$ is infinite. \qed)

In particular, $(F[G/U])_v\otimes_{{\mathcal O}_v}F=F[G/U]$ if and only if 
$U=G$. (Suppose that $U\neq G$ and $F[G/U]$ surjects naturally onto $F[G/U']$ 
for a maximal open subgroup $U'=G_{\{F,L\}|k}$ containing $U$, 
where $L$ is an algebraically closed subfield in $F|k$ of finite 
transcendence degree, and $(F[G/U])_v\otimes_{{\mathcal O}_v}F$ 
surjects naturally onto $(F[G/U'])_v\otimes_{{\mathcal O}_v}F$. 
However, $(F[G/U'])_v\otimes_{{\mathcal O}_v}F=F[\{\text{subfields }F''
\subset{\mathcal O}_v~|~F''\cong L\}]\neq F[G/U']$. \qed) 

\begin{lemma} \label{hens} Let $v:F^{\times}/k^{\times}\longrightarrow
\hspace{-3mm}\to\Gamma$ be a discrete valuation of $F$, and $k\subseteq 
L\subseteq F$ be a subfield. Set $$L^h_v:=\{a\in F~|~\mbox{{\rm for any 
$\gamma\in\Gamma$ there is $a_{\gamma}\in L$ such that} $v(a-a_{\gamma})
>\gamma$}\}\subseteq\overline{L}.$$ Then $W^{G_{F|L^h_v}}=W^{G_{F|L}\cap 
G_v}$ for any smooth $G$-set $W$. However, the inclusion $G_{F|L^h_v}
\supset G_{F|L}\cap G_v$ is strict, unless $\Gamma=0$. \end{lemma}
\begin{proof} Let $w\in W^{G_{F|L}\cap G_v}$. By Proposition \ref{morfizm-mono},
$G_{F|L}\cap G_v\subseteq{\rm Stab}_w\subseteq G_{\{F,K\}|k}$ for $K=\alpha
({\rm Stab}_w)\in A\Pi$. This means that $K\subset\overline{L}$ (otherwise,
for any $f\in K\smallsetminus(K\cap\overline{L})$ and any $t\not\in\overline{L}$ 
with $v(t)\gg 0$ there is $\sigma\in G_{F|\overline{L}}\cap G^{\dagger}_v$ 
such that $\sigma f=f+t\not\in K$), i.e. $G_{F|\overline{L}}
\subseteq G_{F|K}\subseteq{\rm Stab}_w$. On the other hand, $G_{F|L}
\cap G_v\longrightarrow\hspace{-3mm}\to G_{\overline{L}|L^h_v}$, 
and therefore, ${\rm Stab}_w\supseteq G_{F|L^h_v}$. \end{proof} 

\begin{lemma} \label{vychet} In the setup of Lemma \ref{hens} 
the sequence $1\longrightarrow G_{F|L}\cap G_v^{\dagger}
\longrightarrow G_{F|L}\cap G_v\longrightarrow G_{\kappa(v)
|\kappa(v|_L)}\longrightarrow 1$ is exact. \qed \end{lemma}

\begin{lemma} \label{fully-faithfulness} 
If $n=\infty$ then there are canonical isomorphisms 
$${\rm Hom}_G(W,W')\stackrel{\sim}{\longrightarrow}{\rm Hom}_{G_v}(W_v,W')=
{\rm Hom}_{G_v}(W_v,W'_v)$$ for any $W,W'\in{\mathcal S}m_G$. In particular, 
the functor $(-)_v:{\mathcal S}m_G\longrightarrow{\mathcal S}m_{G_v}$, 
is fully faithful. \end{lemma}
\begin{proof} By Lemma \ref{hens}, ${W'}^{G_{F|F'}}={W'}^{G_{F|F'}\cap G_v}$. 
The $G_v$-module $W_v$ is generated by its subspace $W^{G_{F|F'}}$ for any 
algebraically closed subfield $F'\subset{\mathcal O}_v$ of infinite 
transcendence degree, 
so \begin{multline*}{\rm Hom}_{G_v}(W_v,W')\subseteq{\rm Hom}
_{G_{\{F,F'\}|k}\cap G_v}(W^{G_{F|F'}},W')\\ ={\rm Hom}
_{G_{\{F,F'\}|k}\cap G_v}(W^{G_{F|F'}},{W'}^{G_{F|F'}\cap G_v})=
{\rm Hom}_{G_{F'|k}}(W^{G_{F|F'}},{W'}^{G_{F|F'}}).\end{multline*}

This implies that ${\rm Hom}_{G_v}(W_v,W')={\rm Hom}_{G_v}(W_v,W'_v)$. 

As $W\mapsto W^{G_{F|F'}}$ gives an equivalence of categories 
${\mathcal S}m_G\longrightarrow{\mathcal S}m_{G_{F'|k}}$, cf. 
\cite[Lemma 6.7]{repr}, the composition ${\rm Hom}_G(W,W')
\stackrel{\gamma}{\longrightarrow}{\rm Hom}_{G_v}(W_v,W')
\hookrightarrow{\rm Hom}_{G_{F'|k}}(W^{G_{F|F'}},{W'}^{G_{F|F'}})$ 
is an isomorphism, i.e. $\gamma$ is also an isomorphism. \end{proof}

\vspace{4mm}

{\sc Remark.} Clearly, $(-)_v$ does not preserve the irreducibility: the 
surjection $W_v\longrightarrow H_0(G_v^{\dagger},W_v)$ is usually 
non-injective and nontrivial, e.g. $(F/k)_v={\mathcal O}_v/k\supsetneqq
{\mathfrak m}_v=\mathfrak{p}_1\supsetneqq\mathfrak{p}_2\supsetneqq
\dots\supsetneqq\mathfrak{p}_r$ has length $r+1$. However, $(-)_v$ preserves 
the existence of a cyclic vector: if $W\in{\mathcal S}m_G$ is cyclic then 
the $G_{F'|k}$-module $W^{G_{F|F'}}$ admits some cyclic vector $w$ (as 
$H^0(G_{F|F'},-):{\mathcal S}m_G\longrightarrow{\mathcal S}m_{G_{F'|k}}$ 
is an equivalence of categories), and thus, $w$ generates the $G_v$-module 
$W_v$. It follows from Lemma \ref{fully-faithfulness} that if $W$ is 
irreducible and $W_v$ is semisimple then $W_v$ is irreducible. 

\begin{lemma}\label{tochnostt} Let $J$ be a collection of algebraically 
closed subfields $F''\subset F$ of countable transcendence degree over 
$k$. The following conditions on $J$ are equivalent: \begin{itemize}
\item the functor $\Phi_J$ is exact; \item 
for any integer $N\ge 1$, any extension $L$ of $k$ of finite type, 
any collection of embeddings $\xi_j:L\stackrel{/k}{\hookrightarrow}F_j$ 
with $F_j\in J$ for any $1\le j\le N$, and any $\sigma:L
\stackrel{/k}{\hookrightarrow}F$ there is an element 
$\alpha\in{\mathbb Q}[G]$ such that $\alpha\xi_j=0$ for any $1\le j\le N$ 
and $\alpha\sigma-\sigma\in{\mathbb Q}[\{L\stackrel{/k}{\hookrightarrow}
F''~|~F''\in J\}]$; \item for any irreducible $k$-variety $X$, any 
integer $N\ge 1$, any collection of dominant $k$-morphisms $f_j:X
\longrightarrow Y_j$ with $\dim Y_j<\dim X$ for any $0\le j\le N$ 
such that $f_0$ does not factor through $f_j$ for any $1\le j\le N$, 
and any generic point $\sigma:k(Y_0)\stackrel{/k}{\hookrightarrow}F$  
there is a generic 0-cycle $\alpha\in X(F)$ such that $(f_j)_{\ast}
\alpha=0$ for any $1\le j\le N$ and $(f_0)_{\ast}\alpha-\sigma\in
{\mathbb Q}[\{k(Y_0)\stackrel{/k}{\hookrightarrow}F''~|~F''\in J\}]$. 
\end{itemize} \end{lemma}
\begin{proof} $\Phi_J$ is exact if and only if for any exact sequence $0
\longrightarrow A\longrightarrow B\longrightarrow E$ in ${\mathcal S}m_G$ 
its subsequence $\Phi_J(A)\longrightarrow\Phi_J(B)\longrightarrow\Phi_J(E)$ 
is exact, i.e. $\delta_J(A\subset B):=(A\bigcap\Phi_J(B))/\Phi_J(A)$ 
vanishes. We are looking for inclusions $A\subset B$ such that $\delta_J
(A\subset B)\neq 0$. Suppose that an element $\gamma\in A\bigcap\Phi_J(B)$ 
does not belong to $\Phi_J(A)$. Let $\gamma=\sum_{j=1}^N\gamma_j$ 
for some $N\ge 2$ and $\gamma_j\in B^{G_{F|L_j}}$, where $L_j$ are 
some subfields of finite type over $k$ of some elements of $J$. Then 
we can replace $B$ with the submodule generated by $\gamma_1,\dots,
\gamma_N$, and replace $A$ with the submodule generated by $\gamma$. 

Then $A\hookrightarrow B$ is dominated by the diagonal morphism 
$\begin{array}{ccc} C & \stackrel{\Delta}{\longrightarrow} & D:=C^{\oplus N}\\ 
\downarrow&&\phantom{\pi}\downarrow\pi\\ A&\hookrightarrow &B, \end{array}$ 
where $L$ is the compositum of $L_1,\dots,L_N$ and 
$C:={\mathbb Q}[\{L\stackrel{/k}{\hookrightarrow}F\}]$, by sending 
the generator $id:L\stackrel{/k}{\hookrightarrow}F$ of $C$ to $\gamma$ 
and the generator of $j$-th multiple of $D=C^{\oplus N}$ to $\gamma_j$. 

Let $K:=\ker(\pi:D\longrightarrow\hspace{-3mm}
\to B)$. The inclusion $\Phi_J(A)\subset A\bigcap\Phi_J(B)$ is 
dominated by the inclusion $\Phi_J(C)\subset\Delta^{-1}(K+\Phi_J(C)
^{\oplus N})$, since $\Delta^{-1}(K+\Phi_J(C)^{\oplus N})$ coincides 
with $$\{\beta=\beta_j+\varkappa_j\in C~|~\text{$1\le j\le N$, 
$\beta_j\in\Phi_J(C)$, $\varkappa=(\varkappa_1,\dots,\varkappa_N)\in K$}\},$$ 
so $\delta_J(A\subset B)=\Delta^{-1}(K+\Phi_J(D))/
(\Phi_J(C)+\Delta^{-1}(K))$, i.e. $\Phi_J(A)=A\bigcap\Phi_J(B)$ 
if and only if $\Phi_J(C)+\Delta^{-1}(K)=\Delta^{-1}(K+\Phi_J(C)^{\oplus N})$. 

Let $\beta\in\Delta^{-1}(K+\Phi_J(D))\subseteq C$, so 
$\beta=\beta_j+\varkappa_j$ for some $\beta_j\in\Phi_J(C)$, $1\le j\le N$, 
and $\varkappa=(\varkappa_1,\dots,\varkappa_N)\in K$. If $\beta$ has a 
non-zero image in $\delta_J(A\subset B)$ then it has also a non-zero image 
in $\Delta^{-1}(K'+\Phi_J(D))/(\Phi_J(C)+\Delta^{-1}(K'))$, where $K'$ is 
the $G$-submodule in $K$ generated by $\varkappa$. Replacing $A$ with 
$C/\Delta^{-1}(K')$ and $B$ with $D/K'$ we can further suppose that $K$ 
is generated by $\varkappa$. If there is $\alpha\in{\mathbb Q}[G]$ such 
that $\alpha\beta_j=0$ for all $1\le j\le N$ and 
$\alpha\beta-\beta\in\Phi_J(C)$ then $\alpha\varkappa=\Delta\alpha\beta\in 
K$, so $\beta=(\beta-\alpha\beta)+\alpha\beta\in\Phi_J(C)+\Delta^{-1}(K)$.

Conversely, let $A=C/\Delta^{-1}(K)$ and $B=C^{\oplus(N+1)}/K$, where 
$C={\mathbb Q}[\{L\stackrel{/k}{\hookrightarrow}F\}]$ and $K$ generated 
by $\varkappa=(\sigma,\sigma-\xi_1,\dots,\sigma-\xi_N)$. Clearly, 
$$\Delta\sigma=\varkappa+(0,\xi_1,\dots,\xi_N)=\varkappa+
(0,\xi_1,\dots,\xi_N)\in K+\Phi_J(C)^{\oplus(N+1)}.$$ If $\Phi_J$ is 
exact then $\sigma$ belongs, in fact, to $\Delta^{-1}(K)+\Phi_J(C)$, 
i.e. $\sigma=\beta+\psi$ for some $\beta\in\Phi_J(C)$ and 
$\psi\in\Delta^{-1}(K)$. Then, as $\varkappa$ is a generator of $K$, 
there exists $\alpha\in{\mathbb Q}[G]$ such that 
$\alpha\varkappa=(\psi,\dots,\psi)$, and 
therefore, $\sigma-\beta=\alpha\sigma=\alpha(\sigma-\xi_j)$ for any 
$1\le j\le N$, so $\alpha\xi_j=0$ and $\alpha\sigma=\sigma-\beta$. 
Clearly, this is equivalent to the second condition of Lemma. 

The second and the third conditions are related as follows. $X$ is a 
variety with the function field equal to the compositum of all $\xi_j(L)$ 
and $\sigma(L)$, the morphisms $f_j$ are induced by the inclusions of 
the function fields. Clearly, the direct image of 0-cycles to $j$-th multiple 
in the geometric case becomes the action on $\xi_j$ in the algebraic case. 
\end{proof} 

\vspace{4mm}

{\sc `Negative' remarks.} 1. Given a subset $X\subset Y_1\times Y_2$, 
it may well happen that the induced homomorphism $\alpha:{\mathbb Z}[X]
\longrightarrow{\mathbb Z}[Y_1]\times{\mathbb Z}[Y_2]$ is injective, 
even if all the fibres of both projections are infinite. 

To see what is going on, and to construct an example, we associate to 
our data $X\subset Y_1\times Y_2$ a graph $\Gamma$ with coloured edges 
in the following way. The vertices of $\Gamma$ are elements of $Y_1$, 
the colours of the edges of $\Gamma$ are elements of $Y_2$. Any pair 
of vertices is joined by at most one edge of any given colour. A pair 
of vertices $a\neq b$ is joined by an edge of colour $c$ if and only 
if $(a,c),(b,c)\in X$. 

We claim that if $\alpha$ is injective then any pair of vertices is 
joined by at most one edge. Indeed, otherwise, if $a$ and $b$ are joined 
by colours $c$ and $d$ then $(a,c)-(b,c)-(a,d)+(b,d)$ is a non-zero 
element in the kernel of $\alpha$. 

Clearly, if $a$ and $b$, $b$ and $c$ are joined by colour $d$ 
then $a$ and $c$ are joined by colour $d$ (transitivity). 

Suppose that there is a cycle in $\Gamma$, say, with (pairwise distinct) 
vertices $a_i$ and edges of colours $b_i$, $i\in{\mathbb Z}/n{\mathbb Z}$, 
where $a_i$ and $a_{i+1}$ are joined by an edge of colour $b_i$. Then 
$(a_i,b_i),(a_i,b_{i-1})\in X$ for all $i\in{\mathbb Z}/n{\mathbb Z}$, and 
therefore, $\sum_{i\in{\mathbb Z}/n{\mathbb Z}}[(a_i,b_i)-(a_i,b_{i-1})]$ 
is an element in the kernel of $\alpha$. It should be zero, i.e., 
$b_1=\dots=b_n$. In other words, any cycle in $\Gamma$ should be 
monocoloured. 

Similarly, suppose that there is a cycle of colours in $\Gamma$, i.e., 
a collection vertices $a_i$, $1\le i\le n+2$, and (pairwise distinct) 
colours $b_j$, $1\le j\le n$, where $a_i$ and $a_{i+1}$ are joined by an 
edge of colour $b_i$ for $1\le i\le n$ and $a_{n+1}$ and $a_{n+2}$ are 
joined by an edge of colour $b_1$. Then $(a_i,b_i),(a_{i+1},b_i)\in X$ 
for $1\le i\le n$, and $(a_{n+1},b_1),(a_{n+2},b_1)\in X$, so therefore, 
$\sum_{i=2}^n[(a_i,b_i)-(a_{i+1},b_i)]-(a_2,b_1)+(a_{n+1},b_1)$ 
is a non-zero element in the kernel of $\alpha$, leading to contradiction. 
In other words, for any colour $b\in Y_2$ any pair of edges of colour $b$ 
in any connected component of $\Gamma$ has a common vertex. 

Now we produce an example of $\Gamma$. The set of its vertices is 
constructed inductively, as the union of an increasing sequence of sets 
$S_i$. Let $S_0$ be a set (of vertices). We join the elements of $S_0$ 
pairwise by the colour 0. For each $s\in S_0$ fix a set $S_{s,0}$. We 
join the elements of $\{s\}\coprod S_{s,0}$ pairwise by the colour $(s,1)$. 
For each $s\in S_1:=S_0\sqcup\coprod_{s'\in S_0}S_{s,0}$ fix a set 
$S_{s,1}$. We join the elements of $\{s\}\coprod S_{s,1}$ pairwise by 
the colour $(s,2)$. For each $s\in S_2:=S_1\sqcup\coprod_{s'\in S_1}S_{s,1}$ 
fix a set $S_{s,2}$. We join the elements of $\{s\}\coprod S_{s,2}$ pairwise 
by the colour $(s,3)$. And so on. We thus get a graph with vertices 
$Y_1=\bigcup_{j\ge 1}S_j$ and colours of edges 
$Y_2=\{0\}\sqcup\coprod_{j\ge 0}S_j$. 

2. Consider a map $f:X\longrightarrow\prod_{j=1}^NY_j$. It induces 
a homomorphism $\alpha=(f_{j\ast}):{\mathbb Z}[X]\longrightarrow
\prod_{j=1}^N{\mathbb Z}[Y_j]$ such that the diagram 
$$\begin{array}{ccccc}{\mathbb Q}&\stackrel{\deg}{\longleftarrow}
&{\mathbb Q}[X]&\stackrel{id}{\longrightarrow}&{\mathbb Q}[X]\\
\downarrow\lefteqn{\Delta}&&\downarrow\lefteqn{\alpha}&& \| \\
\prod_{j=1}^N{\mathbb Q}&\stackrel{\deg}{\longleftarrow}&
\prod_{j=1}^N{\mathbb Q}[Y_j]&\stackrel{\varphi}{\dashrightarrow}
&{\mathbb Q}[X]\end{array}$$ is commutative. If the induced 
homomorphism $\alpha$ is injective then $f$ is an inclusion 
and there are maps $\varphi_j:Y_j\longrightarrow{\mathbb Q}[X]$ 
such that $\sum_{j=1}^N\varphi_j(f_j(P))=[P]$ for any $P\in X$. 

Let us show that in the situations, we are interested in, one can arrange 
so that there will exist maps $\varphi_j:Y_j\longrightarrow{\mathbb Q}[X]$ 
such that $\sum_{j=1}^N\varphi_j(f_j(P_j))=0$. In other words, we can 
always assume that the induced homomorphism $(f_{j\ast}):{\mathbb Q}[X]^{\circ}
\longrightarrow\prod_{j=1}^N{\mathbb Q}[Y_j]^{\circ}$ is not surjective. 
Changing slightly notations, let an affine variety $X\subset{\mathbb A}^M_k$ 
be given by some polynomial equations, one of which is $P=0$. 
There is an integer $N\ge 2$ and linear functions $f_1,\dots,f_N$ on 
${\mathbb A}^M_k$, and coefficients $a_i\in k$ and $n_i\in{\mathbb N}$ 
such that $P=\sum_ja_if_i^{n_i}$. Choose a ${\mathbb Q}$-linear 
embedding $\lambda:F\hookrightarrow{\mathbb Q}[X(F)]$, and set 
$\varphi_j=\lambda(a_ix^{n_i}):F\longrightarrow{\mathbb Q}[X(F)]$. 

3. The group homomorphisms is a source of non-injectiveness of $\alpha$, cf. 
Proposition \ref{sluchaj-nekotoryh-odnorodnyh-prostranstv} below. One 
could ask, whether in a given situation one can find a group structure on $X$ 
so that the projections to $Y_i$ became homomorphisms. This is not always 
possible. Namely, if $f=(f_1,f_2):X\longrightarrow Y_1\times Y_2$ is a 
homomorphism then $f^{-1}_1(a)=\alpha\ker f_1$, so 
$f_2(f^{-1}_1(a))=\beta f_2(\ker f_1)$. Thus for any $a,b\in Y_1$ one has 
either $f_2(f^{-1}_1(a))=f_2(f^{-1}_1(b))$, or $f_2(f^{-1}_1(a))\cap 
f_2(f^{-1}_1(b))=\emptyset$. In the following example this property does 
not hold, so $f$ is not a homomorphism. Let $L_1,L_2$ be finitely generated 
extensions of $k$, $X=\{L_1\otimes_kL_2\stackrel{/k}{\hookrightarrow}F\}$, 
$Y_j=\{L_j\stackrel{/k}{\hookrightarrow}F\}$ and $f_j:X\longrightarrow Y_j$ 
the restriction maps. Then $f_2(f^{-1}_1(\sigma))=\{\xi:L_2\stackrel{/k}
{\hookrightarrow}F~|~\exists\hat{\xi}:L_1\otimes_kL_2\stackrel{/k}
{\hookrightarrow}F~\text{such that $\hat{\xi}|_{L_1}=\sigma$}\}=\{\xi:L_2
\stackrel{/k}{\hookrightarrow}F~|~\mbox{$\xi(L_2)$ and $\sigma(L_1)$ 
are in general position}\}$, so the intersections of $f_2(f^{-1}_1(\sigma))$ 
and $f_2(f^{-1}_1(\tau))$ are non-empty for all $\sigma,\tau$. On the other 
hand, $f_2(f^{-1}_1(\sigma))\neq f_2(f^{-1}_1(\tau))$ if $\sigma,\tau$ 
are in general position. 

\vspace{4mm}

In view of Lemma \ref{tochnostt} and Proposition 
\ref{sluchaj-nekotoryh-odnorodnyh-prostranstv}, there are some reasons 
to expect that $(-)_v$ is exact. Namely, the second of equivalent conditions 
of Lemma \ref{tochnostt} is satisfied for $\sigma$ in general position with 
respect to the compositum of all $\xi_j(L)$. Also, this condition is 
satisfied in the following situation. Let $L$ be the function field of an 
algebraic $k$-group $A$, and $B$ be an algebraic $k$-subgroup of $A^{N+1}$ 
of dimension $>\dim A$, surjective over all multiples $A$. Suppose that 
$\xi_j$ and $\sigma$ are induced by a generic point $\xi:k(B)\stackrel{/k}
{\hookrightarrow}F$ and the projections $B\subseteq A^{N+1}
\stackrel{{\rm pr}_j}{\longrightarrow}A$. 

({\it Example.} $L$ is purely transcendental of degree one and $\xi_j$ 
sends a generator of $L$ to $x_0+\lambda_j\cdot t$ for some 
$\lambda_j\in k$, $1\le j\le N$, and $\sigma$ sends a generator of $L$ 
to $t$, where $x_0\in{\mathcal O}_v\smallsetminus(k+{\mathfrak m}_v)$ 
and $t\in{\mathfrak m}_v\smallsetminus\{0\}$. Then $A\cong
{\mathbb G}_{a,k}$ and $B\cong{\mathbb G}_{a,k}\times{\mathbb G}_{a,k}$.) 

\begin{proof} We need an element $\alpha\in{\mathbb Q}[G]$ such that 
${{\rm pr}_j}_{\ast}\alpha\xi=0$ for $1\le j\le N$ and 
${{\rm pr}_0}_{\ast}\alpha\xi=0$. 

Assume $\lambda_0=1$ and $\lambda_{N+1}=0$ and that 
$\lambda_j\in k$ are pairwise distinct for $0\le j\le N+1$. Let 
$H\cong({\mathbb Z}/2{\mathbb Z})^{N+1}$ be generated by $e_0,\dots,e_N$ 
and $H_s=H/\langle e_0,\dots,\widehat{e_s},\dots,e_N\rangle$ for all 
$0\le s\le N$. We consider the characters of $H_s$ as characters of $H$: 
$H_s^{\vee}\subset H^{\vee}$. Fix a collection of elements $y_{\chi}\in F$, 
parametrized by the characters $\chi\in H^{\vee}$, linearly independent 
over ${\mathbb Q}^{{\rm ab}}(\lambda_1,\dots,\lambda_N)$ and such that 
\begin{gather}\label{1}y_1+\sum_{s=0}^N\sum_{1\neq\chi\in H_s^{\vee}}
y_{\chi}=x_0\quad\mbox{and}\\ \label{2}\sum_{1\neq\chi\in(H/\langle 
e_1-e_0,\dots,e_N-e_0\rangle)^{\vee}}y_{\chi}-\sum_{s=0}^N
\lambda_s^{-1}\cdot\sum_{1\neq\chi\in H_s^{\vee}}y_{\chi}=t.\end{gather}
For a character $\chi:H\longrightarrow\{\pm 1\}$ set $f_{\chi}=
\sum_{h\in H}\chi(h)\cdot h\in k[H]$, so $gf_{\chi}=\chi(g)\cdot 
f_{\chi}$ for any $g\in H$. 
Let $w=\sum_{h\in H}w_h\cdot h=\sum_{\chi\in(H/\langle e_1-e_0,\dots,
e_N-e_0\rangle)^{\vee}}y_{\chi}\cdot f_{\chi}+\sum_{s=1}^N
(1-\lambda_s^{-1})\cdot\sum_{1\neq\chi\in H_s^{\vee}}y_{\chi}\cdot 
f_{\chi}$ and $x=\sum_{h\in H}x_h\cdot h=y_1\cdot f_1+
\sum_{s=0}^N\sum_{1\neq\chi\in H_s^{\vee}}y_{\chi}\cdot f_{\chi}$. 

One can rewrite the coordinates of $w$ and of $x$ as follows: 
\begin{gather*}w_h=\sum_{\chi\in(H/\langle e_1-e_0,\dots,e_N-e_0
\rangle)^{\vee}}y_{\chi}\cdot\chi(h)+\sum_{s=1}^N(1-\lambda_s^{-1})
\cdot\sum_{1\neq\chi\in H_s^{\vee}}y_{\chi}\cdot\chi(h), \\ 
x_h=y_1+\sum_{s=0}^N\sum_{1\neq\chi\in H_s^{\vee}}y_{\chi}\cdot\chi(h).
\end{gather*}
As $y_{\chi}$ are linearly independent over ${\mathbb Q}^{{\rm ab}}
(\lambda_1,\dots,\lambda_N)$, the coordinates of $w$ and of $x$ 
are pairwise distinct. In particular, their stabilizers 
in ${\mathfrak S}_H$ are trivial. 

For some $1\neq\chi_0\in H^{\vee}_0$ one can choose $y_{\chi}\in F'
\subset{\mathcal O}_v$ algebraically independent over $k(x_0)$ for 
$\chi\neq 1,\chi_0$. Then define $y_{\chi_0}$ by (\ref{2}), and define 
$y_1$ by (\ref{1}). Clearly, $y_1$ is transcendental over 
$k(y_{\chi}~|~\chi\neq 1)$, so (i) $w_h-x_h\in{\mathcal O}_v\smallsetminus
(k+{\mathfrak m}_v)$, (ii) $x_h$ and $w_h$ are algebraically independent 
over $k$ for any $h\in H$, (iii) $x_h$ and $w_{h-e_0}$ are algebraically 
independent over $k$ for any $h\in H$. 

Then we define embeddings $\sigma_h,\tau_h:k(x_0,t)\stackrel{/k}
{\hookrightarrow}F$ by $\sigma_hx_0=\tau_hx_0=x_h$, $\sigma_ht=w_h-x_h$ 
and $\tau_ht=w_{h-e_0}-x_{h-e_0}$ (and extend them arbitrarily to 
elements of $G$). Then the element $\alpha=\sum_{h\in H}(\sigma_h-\tau_h)$ 
satisfies the assumptions of Lemma \ref{tochnostt}. \end{proof}

Slightly more generally, suppose we are given \begin{itemize} 
\item a finite subset $S$ of the projective line ${\mathbb P}^1(k)$, 
\item a finite subset $T\subset{\mathbb Z}
\smallsetminus\{\pm 1\}$ containing a non-zero element, and 
\item a collection of rational functions $(\varphi_j(Y))_{j\in T}$. 
\end{itemize}
Then we get a homomorphism ${\mathbb Z}[F^{\times}\times F\smallsetminus k]
\longrightarrow{\mathbb Z}[F]^{\widetilde{S}}\oplus{\mathbb Z}[F]^T\oplus
{\mathbb Z}[F^{\times}]$, $[(X,Y)]\mapsto\left(([aX-bY])_{(a:b)\in S},
([X^j\varphi_j(Y)])_{j\in T},[X/Y]\right)$, for the pull-back 
$\widetilde{S}$ of $S$ in the set of non-zero pairs in $k$. 

Let us construct a series of elements in its kernel. Let $d$ be the 
least common multiple of the numbers $|j|$ for all non-zero $j\in T$, 
and $P$ be the set of prime divisors of $d$. Replacing $S$ with 
$\bigcup_{\zeta\in\mu_d}\zeta S$, we may assume that $S$ is 
$\mu_d$-invariant. For each line $i\in S$ fix a non-zero vector 
$(a_i,b_i)\in F\times F$ on it. Consider the elements $$A=\sum_{I\subseteq P}
(-1)^{|I|}[(0,\beta)]\stackrel{\oplus}{\cdot}\prod^{\oplus}_{i\in S}\left(
[(0,0)]-[(\prod_{p\in I}\zeta_p\cdot b_i,a_i)]\right),$$ where $\zeta_p$ 
is a primitive $p$-th root of unity and $\stackrel{\oplus}{\cdot}$ and 
$\prod\limits^{\oplus}$ denote the products in the group ring 
${\mathbb Z}[F^2]$ of the additive group $F^2$. 

Set $\zeta_I:=\prod_{p\in I}\zeta_p$. Clearly, the element $A$ is 
equal to $$\sum_{I\subseteq P}\sum_{J\subseteq S}(-1)^{|I|+|J|}
[(\zeta_Ib_J,a_J+\beta)]=\sum_{I\subseteq P}
\sum_{\emptyset\neq J\subseteq S}(-1)^{|I|+|J|}
[(\zeta_Ib_J,a_J+\beta)],$$ where $a_J:=\sum_{i\in J}a_i$ and 
$b_J:=\sum_{i\in J}b_i$, and there are no canceling summands in the 
latter sum for sufficiently general choice of $((a_i,b_i))_{i\in S}$. 
Now, for all appropriate $((a_i,b_i))_{i\in S}$, $\lambda$, $\beta$, 
the element $$\sum_{I\subseteq P}\sum_{\emptyset\neq J\subseteq S}
(-1)^{|I|+|J|}\left([(\zeta_Ib_J,a_J+\beta)]-
[(\lambda\zeta_Ib_J,\lambda a_J+\lambda\beta)]\right)$$ 
belongs to the kernel.

\begin{proposition}\label{sluchaj-nekotoryh-odnorodnyh-prostranstv} Let 
$H$ be an algebraic $k$-group, $N\ge 1$ be an integer, and $H_i$ be a 
$k$-subgroup for each $0\le i\le N$. Suppose that $H_j$ normalizes $H_i$ 
for each pair $1\le i<j\le N$ and $H_i$ is contained in $H_0$ for neither 
$1\le i\le N$. Denote by $f_i:H\longrightarrow H/H_i$ the corresponding 
projections. Then there is a 0-cycle $\alpha\in{\mathbb Q}[H(F)]$ such 
that $(f_i)_{\ast}\alpha=0$ for any $1\le i\le N$ and 
$(f_0)_{\ast}\alpha\neq 0$. More explicitly, almost all 0-cycles of type 
$(h_1-1)\cdots(h_N-1)$, where $h_i\in H_i$ for all $1\le i\le N$, 
satisfy these conditions. \end{proposition}
\begin{proof} The condition $(f_i)_{\ast}\alpha=0$ means that $\alpha\in
\sum_{h\in H_i}{\mathbb Q}[H(F)](h-1)$. It is thus sufficient to check 
that $(h_i-1)\cdots(h_N-1)\in\sum_{h\in H_i}{\mathbb Q}[H(F)](h-1)$. 
This can be shown by descending induction on $i$: 
$(h_i-1)(h_{i+1}-1)=(h_ih_{i+1}h_i^{-1}-1)(h_i-1)+
h_{i+1}(h_{i+1}^{-1}h_ih_{i+1}h_i^{-1}-1)\in\sum_{h\in H_i}
{\mathbb Q}[H(F)](h-1)$, the case $i=N$ being trivial. 

The set $$\{(h_1,\dots,h_N)\in H_1\times\dots\times H_N~|~
(h_1-1)\cdots(h_N-1)\in\sum_{h\in H_0}{\mathbb Q}[H(F)](h-1)\}$$ 
is an algebraic subvariety in $H_1\times\dots\times H_N$. In 
fact, it is a union of subvarieties contained in the subvarieties 
of type $\{(h_1,\dots,h_N)\in H_1\times\dots\times H_N~|~
h_{i_1}\cdots h_{i_s}\in h_{j_1}\cdots h_{j_t}H_0\}$ for all 
$1\le i_1<\dots<i_s\le N$ and some $t\equiv s+1\bmod 2$ and 
$1\le j_1<\dots<j_t\le N$. The latter subvarieties are proper, 
since $H_0$ contains neither of $H_i$. \end{proof}

\begin{lemma} \label{celye-obrazu} If $n=\infty$ and $r=1$ then 
the projection ${\mathbb Q}[\{k(X)\stackrel{/k}{\hookrightarrow}
{\mathcal O}_v\}]\longrightarrow C_{k(X)}$ is surjective 
for any irreducible variety $X$ over $k$. \end{lemma}
\begin{proof} Let $P:k(X)\stackrel{/k}{\hookrightarrow}F$ be 
a generic point. Choose a subfield $k'\subset{\mathcal O}_v
\cap\overline{P(k(X))}$ over $k$ projecting isomorphically 
onto $({\mathcal O}_v\cap\overline{P(k(X))})/
({\mathfrak m}_v\cap\overline{P(k(X))})\subset\kappa(v)$. 

Let $k'(C):=k'P(k(X))\subset F$. If $P$ does not factor through 
${\mathcal O}_v$ then $k'(C)$ is the function field of a smooth 
proper curve over $k'$. As the class of $P$ in $C_{k(X)}$ belongs 
to the image of ${\mathbb Q}[\{k'(C)\stackrel{/k'}{\hookrightarrow}
F\}]\longrightarrow{\mathcal I}_{/k'}{\mathbb Q}[\{k'(C)\stackrel{/k'}
{\hookrightarrow}F\}]={\rm Pic}(C_F)_{{\mathbb Q}}\longrightarrow 
C_{k(X)}$, it remains only to check the surjectivity of ${\mathbb Q}
[\{k'(C)\stackrel{/k'}{\hookrightarrow}{\mathcal O}_v\}]
\longrightarrow{\rm Pic}(C_F)_{{\mathbb Q}}$. 

By \cite[Corollary 3.5]{repr}, the generic points 
$k'(C)\stackrel{/k'}{\hookrightarrow}\kappa(v)$ generate 
${\rm Pic}(C_{\kappa(v)})$. 

Let us check the transitivity of the $(G^{\dagger}_v\cap G_{F|k'})$-action 
on the generic fibres of the specialization ${\rm Pic}(C_F)
\stackrel{s}{\longrightarrow\hspace{-3mm}\to}{\rm Pic}(C_{\kappa(v)})$ 
(i.e. over all $k'({\rm Pic}^jC)\stackrel{\sigma}{\hookrightarrow}
\kappa(v)$). Let $x_1,\dots,x_g$ be a transcendence base of 
$k'({\rm Pic}^jC)$ over $k'$, and $\xi_1,\xi_2:k'({\rm Pic}^jC)
\stackrel{/k}{\hookrightarrow}{\mathcal O}_v$ be a pair of liftings of 
$\sigma$. Then there exists an element $\tau\in G^{\dagger}_v\cap G_{F|k'}$ 
such that $\tau\xi_1 x_j=\xi_2 x_j$ for all $1\le j\le g$. The field 
$k'({\rm Pic}^jC)$ is generated over $k'(x_1,\dots,x_g)$ by an algebraic 
element $f$. According to Hensel's lemma, $\tau\xi_1f=\xi_2f$. (Note, that 
one can replace ${\rm Pic}C$ by an arbitrary smooth variety over $k'$.) 

Thus, the kernel of $s$ coincides with $\{\tau\widetilde{\sigma}
-\widetilde{\sigma}~|~\tau\in G^{\dagger}_v\cap G_{F|k'}\}$, and 
therefore, the generic points $k'(C)\stackrel{/k'}{\hookrightarrow}
{\mathcal O}_v$ (i.e. whose specializations are also generic) generate 
${\rm Pic}(C_F)$. \end{proof}

\begin{corollary} \label{loc-i-glob-invar} $\Gamma(W)=W_v=W$ 
for any $W\in{\mathcal I}_G$. \end{corollary}
\begin{proof} This follows from Lemma \ref{celye-obrazu} in the case 
of $W=C_{k(X)}$. As $C_{k(X)}$ for all irreducible varieties $X$ over 
$k$ form a system of generators of ${\mathcal I}_G$, and the functor 
$W\mapsto W_v$ preserves the surjections, we get that $W_v=W$ for 
arbitrary $W\in{\mathcal I}_G$. \end{proof} 

\begin{lemma} \label{basic-glob} Suppose that $n=\infty$. 
Then $(W_1\otimes W_2)_v\subseteq(W_1)_v\otimes(W_2)_v$ and
$\Gamma(W_1\otimes W_2)\subseteq\Gamma(W_1)\otimes\Gamma(W_2)$ 
for any $W_1,W_2\in{\mathcal S}m_G$. However, $(W\otimes W)_v
\neq W_v\otimes W_v$ if $W={\mathbb Q}[F\smallsetminus k]$. If either $W_1$ is 
a quotient of $A(F)$ for a commutative algebraic $k$-group $A$, 
or $W_1\in{\mathcal I}_G$ then $(W_1\otimes W_2)_v=(W_1)_v
\otimes(W_2)_v$ for any $W_2\in{\mathcal S}m_G$. \end{lemma} 
\begin{proof} The identity 
$(W_1\otimes W_2)^{G_{F|F'}}=W_1^{G_{F|F'}}\otimes W_2^{G_{F|F'}}$ 
(\cite[Lemma 7.5]{pgl}) implies the first inclusion, the second follows 
from the first. If $W={\mathbb Q}[F\smallsetminus k]$, one has 
$(W\otimes W)_v\not\ni[x]\otimes[x']\in W_v\otimes W_v$ for any pair 
of distinct 
$x,x'\in{\mathcal O}_v\smallsetminus(k+{\mathfrak m}_v)$ such that 
$x\equiv x'\bmod{\mathfrak m}_v$. 

Let $w\in W_2^{G_{F|F'}}$. Then $G_{F|F'}\subset G_{F|L}\subseteq
{\rm Stab}_w$, where $L\subset F'$ is of finite type over $k$. 

If $F''$ is an algebraically closed subfield in ${\mathcal O}_v$ 
and either $F''$ and $L$ are algebraically independent over $k$ 
in $\kappa(v)$, or $L\subseteq F''$ then $\overline{LF''}\subset
{\mathcal O}_v$, so $W^{G_{F|F''}}_1\otimes w\subseteq(W_1\otimes W_2)
^{G_{F|\overline{LF''}}}\subseteq(W_1\otimes W_2)_v$. If $W_1$ 
is of type as above then $W^{G_{F|F''}}_1$, with $F''$ and $L$ 
algebraically independent over $k$ in $\kappa(v)$, generate $(W_1)_v$, 
i.e., $(W_1)_v\otimes w\subseteq(W_1\otimes W_2)_v$, thus finally, 
$(W_1)_v\otimes(W_2)_v\subseteq(W_1\otimes W_2)_v$. \end{proof} 

\vspace{4mm}

{\sc Remarks.} 1. If $W$ is endowed with an $F$-vector space 
structure $F\otimes W\longrightarrow W$ then, by Lemma 
\ref{basic-glob}, $W_v$ carries an ${\mathcal O}_v$-module structure: 
$(F\otimes W)_v={\mathcal O}_v\otimes W_v\longrightarrow W_v$. 
Clearly, $F\otimes_{{\mathcal O}_v}W_v\longrightarrow W$ is 
injective, but not surjective, as shows the example of 
$W=F[\{L\stackrel{/k}{\hookrightarrow}F\}]$. 

2. Clearly, $\Gamma_r$ preserves the injections, but not the 
surjections. Namely, let $W:=\bigotimes^N_kF\longrightarrow
\Omega^{N-1}_{F|k}$ be given by $a_1\otimes\dots\otimes a_N
\mapsto a_1da_2\wedge\dots\wedge da_N$. 
Then $W_v=\bigotimes^N_k{\mathcal O}_v$ if $n\ge 2N$, so 
$(\bigotimes^{N-1}_F\Omega^1_{F|k_0})_v=\bigotimes^{N-1}_F
\Omega^1_{{\mathcal O}_v|k_0}$ for any $k_0\subseteq k$; and 
$\Gamma(\bigotimes_k^NF)=k$, but $\Gamma_r(\Omega^{\bullet}_{F|k})=
\Omega^{\bullet}_{F|k,{\rm reg}}$ for any $r\ge 1$. 

\begin{proof} Clearly, $W^{G_{F|F'}}=\bigotimes^N_kF'$, so 
$W_v\subseteq\bigotimes^N_k{\mathcal O}_v$. Let $\omega=
x_1\otimes\dots\otimes x_N$ for some $x_1,\dots,x_N\in F'$ 
whose images in $\kappa(v)$ algebraically independent over 
$k(\overline{y_1},\dots,\overline{y_N})$ for some $y_1,\dots,y_N
\in{\mathcal O}_v$. Then for each $1\le i\le N$ there exists an 
element $\sigma_i\in G_v\cap G_{F|k(x_1,\dots,\widehat{x_i},\dots,
x_N,y_1,\dots,y_N)}$ such that $\sigma_ix_i=x_i+y_i$. Then 
$\prod_{i=1}^N(\sigma_i-1)\omega=y_1\otimes\dots\otimes y_N$. \end{proof} 

In the case $n=\infty$ one can also apply Lemma \ref{basic-glob}. 

\vspace{4mm}

For an integral normal $k$-variety $X$ with $k(X)\subset F$ 
let ${\mathfrak V}(X)$ be the set of all discrete valuations 
of $F$ of rank one trivial on $k$ such that their restrictions 
to $k(X)$ are either trivial, or correspond to divisors on $X$. 
Set ${\mathcal W}(X):=W^{G_{F|k(X)}}\cap
\bigcap_{v\in{\mathfrak V}(X)}W_v\subseteq W$. 

{\sc Remark.} $W^{G_{F|k(X)}}\cap W_v$ depends only on the 
restriction of $v$ to $k(X)$, since the set of $G_{F|k(X)}$-orbits 
$G_{F|k(X)}\backslash G/G_v$ of the valuations of $F$ coincides, 
by Proposition \ref{specializ}, with the set of discrete valuations 
of $k(X)$ of rank $\le r$. E.g., if the restriction of $v$ to 
$k(X)$ is trivial then $W^{G_{F|k(X)}}\subseteq W_v$. 

\vspace{4mm}

{\sc Examples.} 1. If $V={\mathbb Q}[\{L\stackrel{/k}{\hookrightarrow}
F\}]$, or $V=F[\{L\stackrel{/k}{\hookrightarrow}F\}]$ then 
${\mathcal V}(U)=0$ for any non-trivial field extension 
$L|k$ of finite type and any smooth $U$ over $k$. 

2. If $V=\Omega^{\bullet}_{F|k}$ then ${\mathcal V}(U)
=\Omega^{\bullet}_{{\mathcal O}(U)|k}$ for any smooth $U$ over $k$. 

3. If $V={\rm Sym}_F^s\Omega^1_{F|k}$ then ${\mathcal V}(U)\subset
{\rm Sym}_{k(U)}^s\Omega^1_{k(U)|k}$ consists of elements with 
poles (with respect to the lattice ${\rm Sym}_{{\mathcal O}(U)}^s
\Omega^1_{{\mathcal O}(U)|k}$) of order $<s$ for any smooth curve 
$U$ over $k$. 

Note, that ${\mathcal V}$ is functorial with respect 
to all morphisms of smooth $k$-varieties; $\Gamma(V)$ 
is `homotopy invariant' if and only if $s=1$. 

4. If $V=W\otimes F$ for some $W\in{\mathcal I}_G$ 
then ${\mathcal V}(U)=(W^{G_{F|\overline{k(U)}}}\otimes\overline
{{\mathcal O}(U)})^{G_{\overline{k(U)}|k(U)}}$ for any irreducible 
smooth affine $U$ over $k$, where $\overline{{\mathcal O}(U)}$ 
is the integral closure of ${\mathcal O}(U)$ in $F$. 

\vspace{4mm}

Consider the following site $\mathfrak{H}$. Objects of 
$\mathfrak{H}$ are the smooth varieties over $k$. Morphisms in 
$\mathfrak{H}$ are the locally dominant morphisms, transforming 
non-dominant divisors to divisors. Coverings are smooth morphisms 
surjective over the generic point of any divisor on the target. 
Denote by ${\rm Shv}({\mathfrak H})$ the category of sheaves on 
${\mathfrak H}$. Consider the functor $\Phi:{\rm Shv}({\mathfrak H})
\longrightarrow{\mathcal S}m_G$, given by \label{obx-sloj} 
${\mathcal F}\mapsto{\mathcal F}(F):=
\lim\limits_{_A\longrightarrow\phantom{_A}}{\mathcal F}({\bf Spec}(A))
\in{\mathcal S}m_G$. Here $A$ runs over smooth $k$-subalgebras of $F$. 
{\sc Example.} If $j\le 1$ then ${\mathcal F}:
X\mapsto Z^j(X_L)$ is a sheaf on ${\mathfrak H}$, and ${\mathcal F}
(F)=Z^j(L\otimes_kF)$. In particular, $\Phi$ is not faithful, 
since ${\mathcal F}(F)=0$ if $j=1$ and $L=k$. 
\begin{proposition} A choice of $k$-embeddings 
into $F$ of all generic points of all smooth $k$-varieties defines 
a functor ${\mathcal S}m_G\longrightarrow{\rm Shv}({\mathfrak H})$, 
$V\mapsto{\mathcal V}$. \end{proposition} 
{\sc Question.} Is it right adjoint to $\Phi$?

\begin{proof} Clearly, if a dominant morphism $f:U\longrightarrow X$ 
transforms divisors on $U$, non-dominant over $X$, to divisors 
on $X$ then ${\mathfrak V}(U)\subseteq{\mathfrak V}(X)$, 
so ${\mathcal V}(X)\subseteq{\mathcal V}(U)$. 

If, moreover, the pull-back of any divisor on $X$ is 
a divisor on $U$ then ${\mathfrak V}(X)={\mathfrak V}(U)$. 

By Lemma 1.1 of \cite{topologii}, the sequence 
$$0\longrightarrow\prod_{x\in X^0}V^{G_{F|k(x)}}
\longrightarrow\prod_{x\in U^0}V^{G_{F|k(x)}}\longrightarrow
\prod_{x\in(U\times_XU)^0}V^{G_{F|k(x)}}$$ is exact 
(in fact, $X\mapsto\prod_{x\in X^0}V^{G_{F|k(x)}}$ is a sheaf on 
a topology ${\mathfrak D}m_k$). As ${\mathfrak V}(X)={\mathfrak V}(U)=
{\mathfrak V}(U\times_XU)$, the sheaf property for the covering 
$f$ amounts to the exactness of the above sequence restricted 
to $\prod_{x\in U^0}\bigcap_{v\in{\mathfrak V}(X)}V_v$. \end{proof} 

\subsection{The `specialization' functor} \label{specializaciq}
\begin{lemma} \label{dag-coinvar} If $r=1$ and $n=\infty$ 
then $H_0(G^{\dagger}_v,-)=-_{G^{\dagger}_v}$ gives functors 
${\mathcal S}m_G\longrightarrow{\mathcal S}m_{G_{\kappa(v)|k}}$ 
and ${\mathcal I}_G\longrightarrow{\mathcal I}_{G_{\kappa(v)|k}}$. 
There are natural surjections of smooth $G_{\kappa(v)|k}$-modules 
$H^0(G_{F|F'},W)\longrightarrow\hspace{-3mm}\to H_0(G^{\dagger}_v,W_v)$ 
and ${\mathcal I}_{\kappa(v)|k}H_0(G_v^{\dagger},W_v)\longrightarrow
\hspace{-3mm}\to H_0(G^{\dagger}_v,{\mathcal I}W)$ for any $W\in
{\mathcal S}m_G$. They are isomorphisms if $(-)_v$ is exact. \end{lemma} 
\begin{proof} For any smooth representation $W$ of $G$ the stabilizer 
of any vector $\overline{w}\in H_0(G^{\dagger}_v,W)$ contains the 
stabilizer in $G_v$ of its  preimage $w\in W$, which implies the 
smoothness of $H_0(G^{\dagger}_v,W)$, since the projection 
$G_v\longrightarrow\hspace{-3mm}\to G_{\kappa(v)|k}$ is open. 

Choose a presentation of $W$ as cokernel of $\varphi:B\longrightarrow
\bigoplus_{j\in J}{\mathbb Q}[\{L_j\stackrel{/k}{\hookrightarrow}
F\}]$. As the functor $H_0(G^{\dagger}_v,(-)_v)$ is right exact 
if $(-)_v$ is exact, and the functor $H^0(G_{F|F'},-)$ is exact, 
in the following commutative diagram the rows are exact: 
\begin{equation}\label{diagr-specializ}\begin{array}{ccccccc}
H_0(G^{\dagger}_v,B_v)&\longrightarrow & 
\bigoplus_{j\in J}{\mathbb Q}[\{L_j\stackrel{/k}{\hookrightarrow}
\kappa(v)\}]&\longrightarrow&H_0(G^{\dagger}_v,W_v)&\longrightarrow&0\\ 
{\rm onto}\uparrow\phantom{{\rm onto}}&&\phantom{\cong}\uparrow\cong
&&\uparrow\\ B^{G_{F|F'}}&\longrightarrow&\bigoplus_{j\in J}
{\mathbb Q}[\{L_j\stackrel{/k}{\hookrightarrow}F'\}]
&\longrightarrow& W^{G_{F|F'}}&\longrightarrow&0\end{array}\end{equation}
which implies that $ W^{G_{F|F'}}\stackrel{\sim}{\longrightarrow}
H_0(G^{\dagger}_v,W_v)$. 

By \cite[Lemma 6.7]{repr}, the functor $H^0(G_{F|F'},-)$ induces 
equivalences of categories ${\mathcal S}m_G\stackrel{\sim}{\longrightarrow}
{\mathcal S}m_{G_{F'|k}}$ and ${\mathcal I}_G\stackrel{\sim}{\longrightarrow}
{\mathcal I}_{G_{F'|k}}$, so $H_0(G^{\dagger}_v,W_v)\in
{\mathcal S}m_{G_{\kappa(v)|k}}$ and $H_0(G^{\dagger}_v,{\mathcal I}W)
\in{\mathcal I}_{G_{\kappa(v)|k}}$, since they are quotients of 
$W^{G_{F'|k}}\in{\mathcal S}m_{G_{F'|k}}$ and of $({\mathcal I}W)^{G_{F'|k}}
\in{\mathcal I}_{G_{F'|k}}$, respectively. \end{proof} 

\vspace{4mm}

{\sc Remark.} Without assuming exactness of the functor $(-)_v$, it is clear 
from (\ref{diagr-specializ}) that if $\varphi$ is injective then 
$B^{G_{F|F'}}\stackrel{\sim}{\longrightarrow}H_0(G^{\dagger}_v,B_v)$. 
\begin{corollary} For any smooth irreducible divisor $D$ on any smooth proper 
irreducible variety $X$ over $k$ there is a natural morphism $C_{k(D)}
\longrightarrow C_{k(X)}$, if $(-)_v$ is exact for $r=1$, making commutative 
the diagram $\begin{array}{ccc}C_{k(D)}&\longrightarrow\hspace{-3mm}\to & 
CH_0(D_F)_{{\mathbb Q}}\\ \downarrow&&\downarrow\\
C_{k(X)}&\longrightarrow\hspace{-3mm}\to & CH_0(X_F)_{{\mathbb Q}}\end{array}$. 
\end{corollary}
\begin{proof} The second isomorphism from Lemma \ref{dag-coinvar} factors 
through the module ${\mathcal I}_{\kappa(v)|k}H_0(G_v^{\dagger},W)$. Let $W=
{\mathbb Q}[\{k(X)\stackrel{/k}{\hookrightarrow}F\}]$ and ${\mathcal O}_v
\cap k(X)={\mathcal O}_{X,D}$. By Proposition \ref{specializ}, 
$H_0(G_v^{\dagger},W)={\mathbb Q}[\{k(X)
\stackrel{/k}{\hookrightarrow}\kappa(v)\}]\oplus\bigoplus_w{\mathbb Q}
[\{\kappa(w)\stackrel{/k}{\hookrightarrow}\kappa(v)\}]$, where $w$ runs 
over the set of discrete valuations of $k(X)$ of rank 1. Then, replacing 
$\kappa(v)$ by $F$ and applying the functor ${\mathcal I}$, we get a 
morphism $C_{k(X)}\oplus\bigoplus_wC_{\kappa(w)}\longrightarrow C_{k(X)}$. \end{proof} 

\begin{lemma} \label{dag-coinvar-primer} Let ${\mathcal F}$ be a functor 
on the category of smooth $k$-varieties (and all their smooth morphisms). 
For a filtered union ${\mathcal O}=\lim\limits_{_A\longrightarrow\phantom{_A}}
A$ of finitely generated smooth $k$-subalgebras $A$ let 
${\mathcal F}({\mathcal O}):=\lim\limits_{_A\longrightarrow\phantom{_A}}
{\mathcal F}({\bf Spec}(A))$ be the limit. Then \begin{enumerate}
\item ${\mathcal F}({\mathcal O})$ is independent of presentation of 
${\mathcal O}$ as a filtered union; \item ${\mathcal F}(F)^{G_{F|F'}}
={\mathcal F}(F')$ for any $F'=\overline{F'}\subseteq F$ with 
${\rm tr.deg}(F'|k)=\infty$. \end{enumerate} Suppose that ${\mathcal F}
({\mathcal O}_v)={\mathcal F}(F)_v$ for a discrete valuation $v:
F^{\times}/k^{\times}\longrightarrow\hspace{-3mm}\to{\mathbb Q}$.\footnote{The 
valuation ring ${\mathcal O}_v$ of $v$ is the filtered union of the valuation 
rings ${\mathcal O}_L$ of the restrictions of $v$ to subfields $L$ in $F|k$ 
of finite type. As the normal varieties are non-singular outside codimension 
2, the ring ${\mathcal O}_L$ is a filtered union of coordinate rings of 
smooth affine models of $L|k$. 

If $v:F^{\times}/k^{\times}\longrightarrow\hspace{-3mm}\to{\mathbb Q}^r$ is 
a discrete valuation of rank $r>1$ then the valuation ring ${\mathcal O}_v$ 
of $v$ is {\sl not} a union of smooth $k$-subalgebras.} 

Then $H_0(G^{\dagger}_v,{\mathcal F}(F)_v)={\mathcal F}(\kappa(v))$. 
\end{lemma} 
\begin{proof} Fix an isomorphism $\alpha:F\stackrel{\sim}{\longrightarrow}F'$ 
over $k$. Then $\alpha^{\ast}:\lim\limits_{_U\longrightarrow
\phantom{_U}}{\mathcal F}(U)\stackrel{\sim}{\longrightarrow}\lim
\limits_{_V\longrightarrow\phantom{_V}}{\mathcal F}(V)$, 
where ${\mathcal O}(U)\subset F$ and ${\mathcal O}(V)\subset F'$ 
are smooth. For any $U$ there is $\sigma\in G$ such that $\sigma|
_{{\mathcal O}(U)}=\alpha|_{{\mathcal O}(U)}$, so $\lim\limits
_{_V\longrightarrow\phantom{_V}}{\mathcal F}(V)=\lim
\limits_{_{(U,\sigma)}\longrightarrow\phantom{_{(U,\sigma)}}}\sigma
{\mathcal F}(U)={\mathcal F}(F)^{G_{F|F'}}$. 

The composition ${\mathcal F}(F')={\mathcal F}(F)^{G_{F|F'}}
\longrightarrow{\mathcal F}({\mathcal O}_v)\stackrel{i}{\longrightarrow}
{\mathcal F}(\kappa(v))$ is induced by the isomorphism $F'
\stackrel{\sim}{\longrightarrow}\kappa(v)$, so it is itself an isomorphism. 
As $i$ is $G_v^{\dagger}$-invariant, it factors through 
${\mathcal F}({\mathcal O}_v)\longrightarrow\hspace{-3mm}\to 
H_0(G_v^{\dagger},{\mathcal F}({\mathcal O}_v))$, so ${\mathcal F}(F')
\stackrel{\alpha}{\longrightarrow}H_0(G_v^{\dagger},{\mathcal F}
({\mathcal O}_v))$ is injective. Our assumption and 
Lemma \ref{dag-coinvar} imply that $\alpha$ is surjective. \end{proof} 

\vspace{4mm}

{\sc Examples.} 1. As it shows the example of the birationally invariant 
presheaf $U\mapsto\Gamma(\overline{U},{\rm Sym}^2_{{\mathcal O}}
\Omega^1_{\overline{U}|k})$ (or any other stably birational presheaf 
without Galois descent property) and of $F'\subset F$ purely transcendental 
over $k$, one cannot omit the condition that $F'$ is algebraically closed 
to ensure ${\mathcal F}(F)^{G_{F|F'}}={\mathcal F}(F')$. 

2. For a smooth proper $k$-variety $X$ and $q\ge 0$ the functor 
${\mathcal F}:Y\mapsto CH^q(X\times_kY)$ satisfies the assumption of Lemma 
\ref{dag-coinvar-primer} and ${\mathcal F}(F)=CH^q(X_F)$. 

Moreover, let us show that $CH^q(X\times_k{\mathcal O}_v)=CH^q(X_F)$ 
by inducton on $r$. For the minimal prime ideal $0\neq{\mathfrak p}_1
\subset{\mathcal O}_v$, one has a short exact sequence of Gersten complexes of 
flabby sheaves on $X\times_k{\mathcal O}_v$: $$0\to\hspace{-4mm}
\coprod_{x\in(X\times_k({\mathcal O}_v/{\mathfrak p}_1))^{\bullet-1}}
\hspace{-4mm}K_{q-\bullet}(\kappa(x))\to
\coprod_{x\in(X\times_k{\mathcal O}_v)^{\bullet}}
K_{q-\bullet}(\kappa(x))\to\coprod_{x\in X_F^{\bullet}}
K_{q-\bullet}(\kappa(x))\to 0,$$ 
the end of whose long exact cohomological sequence looks as 
$$H^{q-1}(X_F,{\mathcal K}_q)\stackrel{\partial}{\longrightarrow}
CH^{q-1}(X\times_k({\mathcal O}_v/{\mathfrak p}_1))\to 
CH^q(X\times_k{\mathcal O}_v)\stackrel{\pi}{\longrightarrow}CH^q(X_F)\to 0.$$ 
The composition of the surjections $$CH^{q-1}(X_F)\otimes F^{\times}
\stackrel{id\cdot v}{\longrightarrow\hspace{-3mm}\to}CH^{q-1}(X_F)
\stackrel{{\rm sp}}{\longrightarrow\hspace{-3mm}\to}CH^{q-1}(X_L),$$ where 
$L$ is the fraction field of ${\mathcal O}_v/{\mathfrak p}_1$, factors through 
$CH^{q-1}(X_F)\otimes F^{\times}\longrightarrow H^{q-1}(X_F,{\mathcal K}_q)$. 
By induction assumption $CH^{q-1}(X\times_k({\mathcal O}_v/{\mathfrak p}_1))=
CH^{q-1}(X_L)$, which implies the surjectivity of $\partial$, and therefore, 
that $\pi$ is an isomorphism. \qed 

The isomorphism $H_0(G^{\dagger}_v,CH^q(X_F))=CH^q(X_{\kappa(v)})$ is nothing 
but the specialization homomorphism $CH^q(X_F)\longrightarrow\hspace{-3mm}\to 
CH^q(X_{\kappa(v)})$ (cf. \cite{sam}), which is $G^{\dagger}_v$-invariant, 
so it factors through $H_0(G^{\dagger}_v,CH^q(X_F))$. 

2. The functor ${\mathcal F}:Y\mapsto\Gamma(\overline{Y},\Omega^{\bullet}
_{\overline{Y}|k})$, where $\overline{Y}$ is a smooth compactification of $Y$, 
also satisfies the condition of Lemma \ref{dag-coinvar-primer}, and 
${\mathcal F}(F)=\Omega^{\bullet}_{F|k,{\rm reg}}$. 

The reduction modulo the maximal ideal induces a surjection 
$\Omega^{\bullet}_{{\mathcal O}_v|k}\longrightarrow\hspace{-3mm}\to
\Omega^{\bullet}_{\kappa(v)|k}$ and an isomorphism $H_0(G^{\dagger}_v,
\Omega^{\bullet}_{F|k,{\rm reg}})=
\Omega^{\bullet}_{\kappa(v)|k,{\rm reg}}$. \qed 

3. ${\mathcal F}:Y\mapsto{\mathbb Q}[{\mathcal O}(Y)]$ is an example 
with ${\mathcal F}({\mathcal O}_v)={\mathbb Q}[{\mathcal O}_v]\neq
{\mathcal F}(F)_v={\mathbb Q}[{\mathcal O}_v\smallsetminus
(k+{\mathfrak m}_v)]\oplus{\mathbb Q}[k]$. However, still 
$H_0(G^{\dagger}_v,{\mathcal F}(F)_v)={\mathcal F}(\kappa(v))=
{\mathbb Q}[\kappa(v)]$.

\vspace{4mm}

{\sc Remark.} The functor $H_0(G^{\dagger}_v,-):{\mathcal S}m_G
\longrightarrow{\mathcal S}m_{G_{\kappa(v)|k}}$ is neither full, 
since ${\mathbb Q}[\{k(X)\stackrel{/k}{\hookrightarrow}F\}]^G=0$
for $r\ge d>0$, but $H_0(G^{\dagger}_v,{\mathbb Q}[\{k(X)
\stackrel{/k}{\hookrightarrow}F\}])^{G_{\kappa(v)|k}}=
{\mathbb Q}[{\mathcal P}^d_{k(X)}]$, nor faithful, since 
$(F^{\times})^G=k^{\times}$ and $F^G=k$, but 
$H_0(G^{\dagger}_v\cap G_{F|F'},F^{\times})=
H_0(G^{\dagger}_v\cap G_{F|F'},F)=0$.\footnote{since for any $x\in 
F\smallsetminus{\mathcal O}_v$ there are $\sigma,\tau\in G^{\dagger}_v
\cap G_{F|F'}$ such that $\sigma t=tx^{-1}$ and $\tau x=2x$, and 
$F\smallsetminus{\mathcal O}_v$ generates both $F^{\times}$ and $F$.} 
Also one has $H_0(G^{\dagger}_v,V)=0$ for any smooth semilinear 
representation $V$ of $G$ over $F$, if ${\rm tr.deg}(\kappa(v)|k)=\infty$. 

{\it Proof.} Let $w\in V$ and $G_{F|L}\subseteq{\rm Stab}_w$ for some 
$L$ of finite type over $k$. Choose $t\in{\mathfrak m}_v\smallsetminus\{0\}$ 
and $f\in F'$ transcendental over $L(t)$. Then there is 
$\sigma\in G_{F|L(t)}\cap G^{\dagger}_v$ such that $\sigma f-f=t$, 
i.e. $w=\sigma(t^{-1}fw)-t^{-1}fw$. \qed 

\begin{corollary} \label{obr-obr} Let $X$ be an irreducible variety
over $k$ with the function field embedded into $F$ and $Y\subset X$
be an irreducible divisor. The discrete valuations $v:F^{\times}/
k^{\times}\longrightarrow\hspace{-3mm}\to\Gamma$ of $F$ of rank $1$
such that $k(X)\cap{\mathcal O}_v={\mathcal O}_{X,Y}$ (so
$\kappa(v|_{k(X)})=k(Y)$) form a single $G_{F|k(X)}$-orbit.
Then any embedding $k(Y)\stackrel{/k}{\hookrightarrow}F$ induces a
canonical isomorphism $W^{G_{F|k(Y)}}\stackrel{\sim}{\longrightarrow}
H_0(G_v^{\dagger},W_v)^{G_{F|k(X)}\cap G_v}$, if $(-)_v$ is exact. 
\end{corollary}
\begin{proof} By Lemma \ref{dag-coinvar}, $W^{G_{F|F'}}
\stackrel{\sim}{\longrightarrow}H_0(G_v^{\dagger},W_v)$.
Then Lemma \ref{vychet} implies that
$$\left(W^{G_{F|F'}}\right)^{G_{\kappa(v)|\kappa(v|_L)}}\stackrel
{\sim}{\longrightarrow}H_0(G_v^{\dagger},W_v)^{G_{F|L}\cap G_v}
=H_0(G_v^{\dagger},W_v)^{G_{\kappa(v)|\kappa(v|_L)}}.\qquad\qedhere$$
\end{proof}

\subsection{Restrictions on objects of ${\mathcal I}_G$ and 
on the quotients of objects of ${\mathcal I}_G\otimes F$} 
Assume that $n=\infty$. Denote by ${\mathcal C}$ the category of 
smooth semilinear representations of $G$. 
\begin{lemma} The composition ${\mathcal I}_G(k)\stackrel{\otimes_kF}
{\longrightarrow}{\mathcal C}\stackrel{\Gamma\circ{\rm for}}{-\hspace
{-3.5mm}\longrightarrow}{\mathcal S}m_G(k)$ is identical. \end{lemma}
\begin{proof} This follows from Corollary 
\ref{loc-i-glob-invar} and Lemma \ref{basic-glob}. \end{proof} 

\vspace{4mm}

Then one gets the conditions $V_v\otimes_{{\mathcal O}_v}F=V$ 
and $\Gamma(V)\otimes_kF\longrightarrow\hspace{-3mm}\to V$ for any 
$W\in{\mathcal I}_G$ and any semilinear quotient $V$ of $W\otimes F$ 
(`the interesting objects of ${\mathcal C}$ are globally generated'), 
but one has to check, at least, that they are non-empty on the set 
of irreducible objects. Obviously, if the condition 
$V_v\otimes_{{\mathcal O}_v}F=V$ holds for some discrete valuation 
$v$ of rank 1 then it holds for all discrete valuations $v$ of rank 1. 

\begin{corollary} Let ${\mathcal I}'_G$ be the full subcategory of 
${\mathcal I}_G$, for whose objects $W$ the natural map $W^{G_{F|F'}}
\longrightarrow H_0(G^{\dagger}_v,W)$ is an isomorphism, where 
$v\in{\mathcal P}^1_F$ is a discrete valuation of rank 1. If 
$n=\infty$ then ${\mathcal I}'_G$ is an abelian category, 
closed under taking subquotients in ${\mathcal I}_G$. \end{corollary}
\begin{proof} For any $W\in{\mathcal I}'_G$ any any short exact 
sequence $0\longrightarrow W_1\longrightarrow W\longrightarrow W_2
\longrightarrow 0$ the rows of the following commutative diagram 
are exact $$\begin{array}{ccccccccc} 0 & 
\longrightarrow & W^{G_{F|F'}}_1 & \longrightarrow & W^{G_{F|F'}} 
& \longrightarrow & W^{G_{F|F'}}_2 & \longrightarrow & 0 \\ 
&& \alpha_1\downarrow\phantom{\alpha_1} && \phantom{\cong}
\downarrow\cong &&\phantom{\alpha_2}\downarrow\alpha_2 \\ 
&& H_0(G^{\dagger}_v,W_1) & \longrightarrow & H_0(G^{\dagger}_v,W) 
& \longrightarrow & H_0(G^{\dagger}_v,W_2) & \longrightarrow & 0 
\end{array}$$ (the upper one, since $H^0(G_{F|F'},-):{\mathcal S}m_G
\longrightarrow{\mathcal S}m_{G_{F'|k}}$ is an equivalence of categories, 
cf. \cite[Lemma 6.7]{repr}. As $\alpha_1$ and $\alpha_2$ are surjective 
by Lemma \ref{dag-coinvar}, they are isomorphisms. \end{proof} 

\vspace{4mm}

Let ${\mathcal I}^+_G$ (resp., ${\mathcal C}_-$) be the full subcategory 
of ${\mathcal S}m_G$ (resp., of ${\mathcal C}$) with objects $W$ such 
that $W=W_v$ (resp.,  $W=F\otimes_{{\mathcal O}_v}W_v$). Clearly, they 
are closed under taking quotients and contain ${\mathcal I}_G$ (resp., 
${\mathcal I}_G\otimes F$). 
\begin{lemma} Assume that $(-)_v$ is exact. 
Then ${\mathcal I}^+_G$ (resp., ${\mathcal C}_-$) is a Serre 
subcategory in ${\mathcal S}m_G$ (resp., in ${\mathcal C}$). 
Moreover, ${\mathcal I}^+_G\neq{\mathcal I}_G$. 

The inclusion functors ${\mathcal I}^+_G\hookrightarrow
{\mathcal S}m_G$ and ${\mathcal C}_-\hookrightarrow{\mathcal C}$ 
admit right adjoints $W\mapsto\Gamma(W)$ and $V\mapsto\bigcap_v
(F\otimes_{{\mathcal O}_v}V_v)$, respectively, but do not admit 
left adjoints. \end{lemma}
\begin{proof} The first assertion is clear, since the functors $(-)_v$ 
and $F\otimes_{{\mathcal O}_v}(-)_v$ are exact. 

Let $W'\in{\mathcal I}^+_G$, $\varphi\in{\rm Hom}_{{\mathcal S}m_G}(W',W)$ 
and $V'\in{\mathcal C}_-$, $\psi\in{\rm Hom}_{{\mathcal C}}(V',V)$. Then 
$\varphi$ factors through $\Gamma(W)$ and $\psi$ factors through 
$\bigcap_v(F\otimes_{{\mathcal O}_v}V_v)$. As $(-)_v$ is exact, 
$\Gamma(W)_v=\Gamma(W)\cap W_v=\Gamma(W)$ and $F\otimes_{{\mathcal O}_u}
(\bigcap_v(F\otimes_{{\mathcal O}_v}V_v))_u=\bigcap_v
(F\otimes_{{\mathcal O}_v}V_v)\cap(F\otimes_{{\mathcal O}_u}V_u)
=\bigcap_v(F\otimes_{{\mathcal O}_v}V_v)$. 
Moreover, $\Gamma({\rm Sym}^2_F\Omega^1_{F|k})_v=
\Gamma({\rm Sym}^2_F\Omega^1_{F|k})\not\in{\mathcal I}_G$. 

To see the absence of left adjoints, let us check that ${\mathbb Q}
[F\smallsetminus k]\not\in{\mathcal I}^+_G$ can be embedded into a product 
of copies of $\prod_{j\ge 2}\Gamma({\rm Sym}^j_F\Omega^1_{F|k})$. For each 
integer $s\ge 2$ consider ${\mathbb Q}[F\smallsetminus k]\longrightarrow
\prod_{j\ge 2}{\rm Sym}^j_F\Omega^1_{F|k}$ given by $[x]\mapsto
\left(\frac{(dx)^j}{x^{j-1}(x^s-1)}\right)_{j\ge 2}$. Suppose that an 
element $\sum_ia_i[x_i]$ with pairwise distinct $x_i\in F\smallsetminus k$ 
is sent to zero for any $s$. Then, by Artin's theorem on linear 
independence of characters, $\sum_{i:~x_i/x\in k^{\times}}
a_i\frac{x_i}{x_i^s-1}=0$ for any 
$x\in F^{\times}$. Fix some $i$. Then $\sum_{j:~x_j=\lambda_jx_i}a_i
\frac{\lambda_j}{\lambda_j^sx_i^s-1}=0$. As $\lambda_j\in k^{\times}$ 
are pairwise distinct, for $s=1+N!$ and sufficiently big $N$ the left 
hand side has a simple pole at $x_i=1$, unless $a_i=0$. \end{proof} 

\vspace{4mm}

{\sc Remarks.} 1. Assuming that Corollary \ref{obr-obr} holds, the following 
construction should provide a fully faithful functor from ${\mathcal I}^+_G$ 
to the category of (birationally invariant) functors on smooth varieties 
over $k$ with all, not necessarily smooth, morphisms, which is right 
quasi-inverse to the functor $\Phi:{\mathcal F}\mapsto{\mathcal F}(F)$, 
cf. p.\pageref{obx-sloj}. 

As usually, we assume that the function fields of irreducible 
varieties $Y\subset X$ over $k$ are embedded into $F$. For any $W\in
{\mathcal I}^+_G$ the natural homomorphism $W^{G_{F|k(X)}}\longrightarrow 
H_0(G^{\dagger}_v,W)$ factors through $W^{G_{F|k(X)}}\longrightarrow 
H_0(G^{\dagger}_v,W_v)^{G_v\cap G_{F|k(X)}}$. By Corollary \ref{obr-obr}, 
the target space is canonically isomorphic to $W^{G_{F|k(Y)}}$, if 
$k(X)\cap{\mathcal O}_v={\mathcal O}_{X,Y}$. \qed 

2. Let $W\in{\mathcal C}$ be a sub-object of $F[\{L\stackrel{/k}
{\hookrightarrow}F\}]$ for a finitely generated subfield $L$ in $F|k$. 
The following conditions are equivalent. \begin{enumerate} \item 
$F[\{L\stackrel{/k}{\hookrightarrow}F\}]/W\in{\mathcal C}_-$; \item $W
+F[\{L\stackrel{/k}{\hookrightarrow}{\mathcal O}_v\}]=F[\{L\stackrel
{/k}{\hookrightarrow}F\}]$; \item $W\longrightarrow\hspace{-3mm}\to 
F[\{L\stackrel{/k}{\hookrightarrow}F\}]/F[\{L\stackrel{/k}
{\hookrightarrow}{\mathcal O}_v\}]$; \item for any $\sigma:L\stackrel
{/k}{\hookrightarrow}F$ there exists an element $\alpha_{\sigma}\in 
F[\{L\stackrel{/k}{\hookrightarrow}{\mathcal O}_v\}]$ such that 
$\sigma+\alpha_{\sigma}\in W$. \end{enumerate} Now, given any 
$G_v$-equivariant map $\alpha:\{L\stackrel{/k}{\hookrightarrow}F\}
\smallsetminus\{L\stackrel{/k}{\hookrightarrow}{\mathcal O}_v\}
\longrightarrow F[\{L\stackrel{/k}{\hookrightarrow}{\mathcal O}_v\}]$, 
we define $W_{\alpha}\in{\mathcal C}$ as the sub-object generated by 
$\sigma+\alpha_{\sigma}$ for all $\sigma$ in the domain of definition of 
$\alpha$. The objects $F[\{L\stackrel{/k}{\hookrightarrow}F\}]/W_{\alpha}$ 
for all $\alpha$ and purely transcendental extensions $L$ of $k$ of all 
(sufficiently big) finite transcendence degrees form a system of 
generators of the category ${\mathcal C}_-$. 

\vspace{5mm}

\noindent
{\sl Acknowledgement.} {\small I would like to thank Uwe Jannsen 
for many inspiring discussions.  
I am grateful to the referee for pointing out several inaccuracies 
in the previous version of the paper and to Viktor Kulikov for 
explaining to me that the valuation ring of a discrete valuation $v$ is 
a union of smooth $k$-subalgebras only if rank of $v$ is equal to one. 
I am grateful to Regensburg University for hospitality, and to 
Alexander von Humboldt-Stiftung for support that made my stay in 
Regensburg possible. 

}

\end{document}